\documentclass{emsprocart}

\usepackage{srcltx}
\usepackage{stmaryrd}
\usepackage[all]{xy}

\contact[leclerc@math.unicaen.fr]{Bernard Leclerc,\\ LMNO UMR 6139 CNRS, Universit\'e de Caen,\\ Institut Universitaire de France,\\ F-14032 Caen cedex, France}

\newtheorem{theorem}{Theorem}[section]
\newtheorem{corollary}[theorem]{Corollary}

\newtheorem{proposition}[theorem]{Proposition}
\newtheorem{conjecture}[theorem]{Conjecture}

\theoremstyle{definition}

\newcommand{\Hom}{\operatorname{Hom}}
\newcommand{\Ker}{\operatorname{Ker}}
\newcommand{\Id}{\operatorname{Id}}
\newcommand{\im}{\operatorname{Im}}
\newcommand{\End}{\operatorname{End}}
\newcommand{\Gr}{\operatorname{Gr}}
\newcommand{\FM}{\operatorname{FM}}
\newcommand{\bd}{\mathbf{d}} 
\def\N{{\mathbb N}}
\def\Z{{\mathbb Z}}
\def\C{{\mathbb C}}
\def\Q{{\mathbb Q}}
\def\P{{\mathbb P}}
\def\CC{\mathcal{C}}
\def\Y{\mathcal{Y}}
\def\M{\mathcal{M}}
\def\A{\mathcal{A}}
\def\F{\mathcal{F}}
\def\D{\mathcal{D}}
\def\NN{\mathcal{N}}
\def\gg{\mathfrak g}
\def\eg{\emph{e.g.\,}}
\def\ie{\emph{i.e.\,}}
\def\resp{\emph{resp.\,}}
\def\Sl{\mathfrak{sl}}
\def\so{\mathfrak{so}}
\def\hSl{\widehat{\mathfrak{sl}}}
\def\hg{\widehat{\mathfrak g}}
\def\hI{\widehat{I}}
\def\hL{\widehat{\Lambda}}
\def\al{\alpha}
\def\be{\beta}
\def\ga{\gamma}

\def\vpi{\varpi}
\def\om{\omega}
\def\la{\lambda}
\def\Si{\Sigma}
\def\L{\Lambda}
\def\fM{\mathfrak{M}}
\def\fL{\mathfrak{L}}
\def\tF{\widetilde{\F}}
\def\De{\Delta}
\def\Ga{\Gamma}
\def\<{\langle\,}
\def\>{\,\rangle}
\def\x{\mathbf{x}}
\def\y{\mathbf{y}}
\def\z{\mathbf{z}}
\def\ds{\displaystyle}
\def\tQ{\widetilde{Q}}
\def\si{\sigma}

\title[Quantum loop algebras, quiver varieties, and cluster algebras]{Quantum loop algebras, quiver varieties, and \\
cluster algebras}

\author[Bernard Leclerc]{Bernard Leclerc}

\begin{document}

\begin{abstract}
These notes reflect the contents of three lectures given at the
workshop of the 14th International Conference on Representations
of Algebras (ICRA XIV), held in August 2010 in Tokyo. 
We first provide an introduction to quantum 
loop algebras and their finite-dimensional representations.
We explain in particular Nakajima's geometric description of the irreducible
$q$-characters in terms of graded quiver varieties.
We then present a recent attempt to understand the tensor
structure of the category of finite-dimensional representations
by means of cluster algebras. This takes the form of a general
conjecture depending on a level $\ell\in\N$. The conjecture
for $\ell = 1$ is now proved thanks to some joint work with Hernandez,
and a subsequent paper of Nakajima. The general case is still open.
\end{abstract}



\maketitle

\section*{Introduction}

At the very origin of the theory of quantum groups is the search of 
an algebraic procedure for constructing solutions of the quantum
Yang-Baxter equation
\begin{equation}\label{YB}
 R_{12}(u)R_{13}(uv)R_{23}(v) =  R_{23}(v)R_{13}(uv)R_{12}(u).
\end{equation}
The unknown $R(u)$ of this equation
is an endomorphism of $V\otimes V$ for
some finite-dimensional vector space $V$.
This endomorphism depends on a
parameter $u\in\C^*$,
and we denote by $R_{ij}(u)$
the endomorphism of $V\otimes V\otimes V$ acting via
$R(u)$ on the product of the $i$th and $j$th factors, and via $\Id_V$ 
on the remaining factor. 

The Yang-Baxter equation appeared in several different guises in the literature
on integrable systems (see \cite{J1} for a nice selection of early
papers on the subject).
As is well known, Drinfeld and Jimbo showed how to associate a solution of
(\ref{YB}) with every irreducible finite-dimensional representation $V$ of a quantum
loop algebra $U_q(L\gg)$. Here $\gg$ is a simple complex Lie algebra, 
$L\gg = \gg\otimes\C[t,t^{-1}]$ is its loop algebra, and $U_q(L\gg)$ denotes
the quantum analogue of the enveloping algebra of $L\gg$ with parameter
$q\in\C^*$ not a root of unity. 

This gives a strong motivation for studying the category $\CC$ of finite-dimensional 
representations of $U_q(L\gg)$, and many authors have brought important contributions
(see \eg \cite{AK,CP2,CP,FR,GV,N1}).
In particular the simple objects of $\CC$ have been classified by 
Chari and Pressley, an appropriate notion of $q$-character has been
introduced by Frenkel and Reshetikhin, and, when $\gg$ is of simply laced type,
Nakajima has calculated the irreducible $q$-characters 
in terms of the cohomology of certain quiver varieties.

After reviewing these results in the first two lectures,
we will turn to more recent attempts
to understand the tensor structure of $\CC$. In the case of 
$\gg = \Sl_2$, Chari and Pressley \cite{CP2} proved that every
simple object is a tensor product of Kirillov-Reshetikin modules.
These irreducible representations (which are defined for every $\gg$)
are very particular. 
They also come from the theory of integrable systems
\cite{KR, KNS} where they were first studied in relation with row-to-row transfer matrices.
However, if $\gg$ is different from $\Sl_2$, the Kirillov-Reshetikin modules are not
the only prime simple objects and the situation is far more complicated.
Thus, already for $\gg=\Sl_3$ we do not know a factorization 
theorem for simple objects. 

In the last lecture, we will explain a general conjecture
(for $\gg$ of simply-laced type) which would imply that
the prime tensor factorization of many simple objects can be described
as the factorization of a cluster monomial into a product of cluster variables. 
In retrospect, this should not be too surprising, since one of the
first applications of cluster combinatorics given by Fomin and Zelevinsky
was the proof of Zamolodchikov's periodicity conjecture for $Y$-systems \cite{FZ1}, 
also intimately related with the representation theory of $U_q(L\gg)$ \cite{KNS}.
Our conjecture, which involves a positive integer $\ell$, is now proved
when $\ell = 1$ \cite{HL, N5}, but the general case remains open. 

\section{Representations of quantum loop algebras}

In this first lecture, we review the definition of quantum loop algebras,
the classification of their finite-dimensional irreducible representations,
and the notion of $q$-character. 
For simplicity, we only consider quantum loop algebras of simply laced
type $A, D, E$.
We also formulate the $T$-systems satisfied by the $q$-characters of
the Kirillov-Reshetikhin modules. All this is illustrated in the case
of $U_q(L\Sl_2)$. 
We briefly explain the connection with the Yang-Baxter equation.
Finally, we introduce some interesting tensor subcategories of the category 
of finite-dimensional representations, which will be used in \S\ref{sect_4}. 

For a recent and more complete survey on these topics, see \cite{CH}.

\subsection{The quantum loop algebra $U_q(L\gg)$} \label{subsect_2.1}

Let $\gg$ be a simple Lie algebra over $\C$ of type $A_n, D_n$ or $E_n$.
We denote by $I=[1,n]$ the set of vertices of the Dynkin diagram,
by $A=[a_{ij}]_{i,j\in I}$ the Cartan matrix, and
by $\Pi=\{\al_i\mid i\in I\}$ the set of simple roots.
Let $L\gg = \gg\otimes\C[t,t^{-1}]$ be the loop algebra of $\gg$.
This is a Lie algebra with bracket
\begin{equation}
[x\otimes t^k,\, y\otimes t^l] = [x,y] \otimes t^{k+l}, \quad (x,y\in\gg,\ k,l\in \Z). 
\end{equation}
Following Drinfeld \cite{D}, the enveloping algebra $U(L\gg)$ has 
a quantum deformation $U_q(L\gg)$.
This is an algebra over $\C$ defined by a presentation 
with infinitely many generators 
\begin{equation}
x^+_{i,r},\ x^-_{i,r},\ h_{i,m},\ k_i,\ k_i^{-1}, 
\qquad (i\in I,\ r\in\Z,\ m\in\Z\setminus \{0\}), 
\end{equation}
and a list of relations which we will not repeat (see \eg \cite{FR}).
These relations depend on $q\in\C^*$ which we assume is not a
root of unity.
If $x^+_i,\,x^-_i,\,h_i\ (i\in I)$ denote the Chevalley generators
of $\gg$ then $x^\pm_{i,r}$ is a $q$-analogue of $x^\pm_i\otimes t^r$,
$h_{i,m}$ is a $q$-analogue of $h_i\otimes t^m$, and
$k_i$ stands for the $q$-exponential of $h_i\equiv h_i\otimes 1$.
In fact $U_q(L\gg)$ is isomorphic to a quotient of the quantum
enveloping algebra $U_q(\hg)$ attached by Drinfeld and Jimbo 
to the affine Kac-Moody algebra $\hg$.
It thus inherits from $U_q(\hg)$ the structure of a Hopf algebra. 

Our main object of study is the category $\CC$ of finite-dimensional 
$U_q(L\gg)$-modu\-les\footnote{We only consider modules \emph{of type 1}, a mild technical
condition, see \eg \cite[\S12.2 B]{CP}.}.
Since $U_q(L\gg)$ is a Hopf algebra, $\CC$ is an abelian monoidal category.
It is well-known that $\CC$ is not semisimple.
We denote by $R$ its Grothendieck ring.
Given two objects $M$ and $N$ of $\CC$, the tensor products $M\otimes N$
and $N\otimes M$ are in general not isomorphic. 
However, they have the same composition factors with the same multiplicities,
so $R$ is a commutative ring. 

For every $a\in\C^*$ there exists an automorphism $\tau_a$ of $U_q(L\gg)$
given by
\[
\tau_a(x^\pm_{i,r}) = a^r x^\pm_{i,r},\quad
\tau_a(h_{i,m}) = a^m h_{i,m},\quad
\tau_a(k_i^{\pm1}) = k_i^{\pm1}.
\]
(This is a quantum analogue of the automorphism $x\otimes t^k \mapsto a^k(x\otimes t^k)$
of $L\gg$.)
Each automorphism $\tau_a$ induces an auto-equivalence $\tau_a^*$ of $\CC$, which maps  
an object $M$ to its pullback $M(a)$ under $\tau_a$.

\subsection{$q$-characters}
By Drinfeld's presentation, the generators $k_i^{\pm1}$ and $h_{i,m}$
are pairwise commutative.  
So every object $M$ of $\CC$ can be written as a finite direct sum of common 
generalized eigenspaces for the simultaneous action of the $k_i$ and of the $h_{i,m}$.
These common generalized eigenspaces are called the
{\em l}-weight-spaces of $M$. (Here {\em l} stands for ``loop'').
The $q$-character of $M$, introduced by Frenkel and Reshetikhin \cite{KR}, 
is a Laurent polynomial with positive
integer coefficients in some indeterminates
$Y_{i,a}\ (i\in I, a\in\C^*)$,
which encodes the decomposition of~$M$
as the direct sum of its {\em l}-weight-spaces.
More precisely, the eigenvalues of the $h_{i,m}\ (m>0)$ 
in an {\em l}-weight-space $V$ of $M$ 
are always of the form
\begin{equation}\label{eq24}
\frac{q^m-q^{-m}}{m(q-q^{-1})}
\left(
\sum_{r=1}^{k_i}(a_{i,r})^m-\sum_{s=1}^{l_i}(b_{i,s})^m
\right)
\end{equation}
for some nonzero complex numbers $a_{i,r}, b_{i,s}$.
Moreover, they completely determine 
the eigenvalues of the $h_{i,m}\ (m<0)$ and of the $k_i$ on $V$.
We encode this collection of eigenvalues with 
the Laurent monomial
\begin{equation}\label{eq23}
m_V = \prod_{i\in I}\left(
\prod_{r=1}^{k_i} Y_{i,a_{i,r}} \prod_{s=1}^{l_i} Y_{i,b_{i,s}}^{-1}
\right).
\end{equation}
The collection of eigenvalues (\ref{eq24}), or equivalently 
the monomial (\ref{eq23}),
will be called the {\em l}-weight
of $V$. 
Let $\Y = \Z[Y_{i,a}^{\pm1} ; i\in I, a\in\C^*]$. 
One then defines 
the $q$-character of $M\in\CC$ by
\begin{equation}
\chi_q(M) = \sum_V \dim V \, m_V \in \Y,
\end{equation}
where the sum is over all {\em l}-weight spaces $V$ of $M$
(see \cite[Prop. 2.4]{FM}).
\begin{theorem}[\cite{FR}]
The Laurent polynomial $\chi_q(M)$ depends only on the class of $M$
in $R$, and the induced map
$\chi_q : R \to \Y$ is an injective ring homomorphism. 
\end{theorem}
The subalgebra of $U_q(L\gg)$ generated by 
\[
x^+_{i,0},\ x^-_{i,0},\ \ k_i,\ k_i^{-1}, 
\qquad (i\in I), 
\]
is isomorphic to $U_q(\gg)$.
Hence every $M\in\CC$ can be regarded as a $U_q(\gg)$-module by restriction.
The {\em l}-weight-space decomposition of $M$ is a refinement
of its decomposition as a direct sum of $U_q(\gg)$-weight-spaces. 
Let $P$ be the weight lattice of $\gg$, with basis given by the fundamental
weights $\vpi_i\ (i\in I)$.
Let $\Z[P]$ be the group ring of $P$. 
As usual, $\la\in P$ is written in $\Z[P]$
as a formal exponential $e^\la$ to allow multiplicative notation.
We denote by $\om$ the ring homomorphism from $\Y$ to $\Z[P]$ defined by
\begin{equation}
\om\left(Y_{i,a}\right) = e^{\vpi_i}.
\end{equation}
If $V$ is an {\em l}-weight-space of $M$ with {\em l}-weight the 
Laurent monomial $m\in \Y$, then $V$ is a subspace of the $U_q(\gg)$-weight-space
with weight $\la$ such that $e^\la=\om(m)$.
Hence, the image $\om(\chi_q(M))$ of the $q$-character of $M$ is
the ordinary character of the underlying $U_q(\gg)$-module.

For $i\in I$ and $a\in\C^*$ define
\begin{equation}\label{eqrootA}
 A_{i,a} = Y_{i,aq}Y_{i,aq^{-1}}\prod_{j\not = i}Y_{j,a}^{a_{ij}}.
\end{equation}
Thus $\om(A_{i,a})=e^{\al_i}$, and the $A_{i,a}\ (a\in \C^*)$ should be viewed as affine 
analogues of the simple root~$\al_i$\footnote{In the non-simply-laced case, the definition
of $A_{i,a}$ is more complicated, see \cite{FR}.}.
Following \cite{FR}, we define a partial order on the set $\M$ 
of Laurent monomials in the variables $Y_{i,a}$ by setting:
\begin{equation}\label{eqorder}
m \le m'\quad \Longleftrightarrow \quad 
\mbox{$\displaystyle\frac{m'}{m}$ is a monomial in the $A_{i,a}$ with exponents $\ge 0$.}
\end{equation}
This is an affine analogue of the usual partial order on $P$,
defined by $\la \le \la'$ if and only if $\la'-\la$ is a sum of simple roots $\al_i$.

A monomial $m\in\M$ is called dominant if it does not contain
negative powers of the variables $Y_{i,a}$. 
We will denote by $\M_+$ the set of dominant monomials.
They parametrize simple objects of $\CC$, as was first shown
by Chari and Pressley \cite{CP}\footnote{The original parametrization 
of \cite{CP} is in terms of Drinfeld polynomials, but in these notes 
we will rather use the equivalent parametrization
by dominant monomials.}. More precisely, we have
\begin{theorem}[\cite{FM}]
Let $S$ be a simple object of $\CC$. 
The $q$-character of $S$ is of the form
\begin{equation}\label{eq_q_char}
 \chi_q(S) = m_S\left(1 + \sum_p M_p\right),
\end{equation}
where $m_S\in\M_+$, and all the $M_p\not = 1$ are monomials in the variables $A_{i,a}^{-1}$
with nonnegative exponents. Moreover the map $S \mapsto m_S$ induces a bijection
from the set of isoclasses of irreducible modules in $\CC$ to $\M_+$.
\end{theorem}
The dominant monomial $m_S$ is called
the highest {\em l}-weight
of $\chi_q(S)$, since every other monomial $m_SM_p$ of (\ref{eq_q_char}) is
less than $m_S$ in the partial order (\ref{eqorder}).
The one-dimensional {\em l}-weight-space of $S$ with {\em l}-weight
$m_S$ consists of the highest-weight vectors of $S$, that is,
the {\em l}-weight vectors $v\in S$ such that $x_{i,r}^+v=0$ for every 
$i\in I$ and $r\in\Z$. 

For $m\in\M_+$, we denote by $L(m)$ the corresponding simple object
of $\CC$. In particular, the modules 
\[
L(Y_{i,a}),\qquad (i\in I,\ a\in \C^*),
\]
are called the fundamental modules.
It is known \cite[Cor. 2]{FR} that the Grothendieck ring $R$ is the polynomial 
ring over $\Z$ in the classes of the fundamental modules.

\subsection{Kirillov-Reshetikhin modules}\label{subsect_2.3}

For $i\in I$, $k\in \N^*$ and $a\in\C^*$, the simple object 
$W^{(i)}_{k,a}$ with highest {\em l}-weight
\[
m_{k,a}^{(i)} = \prod_{j=0}^{k-1}Y_{i,\,aq^{2j}}
\]
is called a Kirillov-Reshetikhin module.
Thus $W^{(i)}_{k,a}$ is an affine analogue of the irreducible
representation of $U_q(\gg)$ with highest weight $k\vpi_i$.
In particular for $k=1$,
$W^{(i)}_{1,a}$ coincides with the fundamental module $L(Y_{i,a})$.
By convention, $W^{(i)}_{0,a}$ is the trivial representation for every
$i$ and $a$.

The classes $[W^{(i)}_{k,a}]$ in $R$, or equivalently the 
$q$-characters $\chi_q(W^{(i)}_{k,a})$,
satisfy the following system of equations
indexed by $i\in I$, $k\in\N^*$, and $a\in\C^*$, called the $T$-system\footnote{In 
the non-simply laced case, the $T$-systems are more complicated, see \cite{KNS,H2}.}: 
\begin{equation}\label{eqTsystem}
[W^{(i)}_{k,a}][W^{(i)}_{k,aq^2}] = [W^{(i)}_{k+1,a}][W^{(i)}_{k-1,aq^2}]
+ \prod_{j\not = i} [W^{(j)}_{k,aq}]^{-a_{ij}}. 
\end{equation}
This was conjectured in \cite{KNS} and proved in \cite{N2, H2}.
Using these equations, one can calculate inductively the expression
of any $[W^{(i)}_{k,a}]$ as a polynomial in the classes 
$[W^{(i)}_{1,a}]$ of the fundamental modules.
Thus, one can obtain the $q$-characters of all the Kirillov-Reshetikhin
modules once the $q$-characters of the fundamental modules are known.

\subsection{The case of $\Sl_2$}\label{subsect_2.4}

Let us illustrate the previous statements for $\gg=\Sl_2$. 
Here $I=\{1\}$, and we may drop the index~$i$ in 
$\al_i$, $\vpi_i$, $Y_{i,a}$, $A_{i,a}$, $[W^{(i)}_{k,a}]$.
We have $A_a = Y_{aq}Y_{aq^{-1}}$.
The fundamental modules $W_{1,a}\ (a\in\C^*)$ are the
affine analogues of the vector representation 
$\C^2$ of $U_q(\Sl_2)$, whose character is equal to 
\[
e^{\vpi}+e^{-\vpi} = e^\vpi(1+e^{-\al}).
\] 
Since the {\em l}-weight spaces are subspaces of the one-dimensional
$U_q(\Sl_2)$-weight spaces,
$W_{1,a}$ also decomposes as a sum of 
two {\em l}-weight spaces, and the two {\em l}-weights are easily
checked to be $Y_{a}$ and
$Y_{aq^2}^{-1}$. Hence
\[ 
\chi_q(W_{1,a})=Y_a+Y_{aq^2}^{-1}=Y_a(1+A_{aq}^{-1}).
\]
The $T$-system (\ref{eqTsystem}) reads
\[
[W_{k,a}][W_{k,aq^2}] = [W_{k+1,a}][W_{k-1,aq^2}]
+ 1, \qquad (a\in\C^*,\ k\in\N^*). 
\]
From the identity 
\begin{equation}\label{eqKR1}
[W_{1,a}][W_{1,aq^2}] = [W_{2,a}][W_{0,aq^2}]+ 1=[W_{2,a}]+1, 
\end{equation}
one deduces that
\[
\chi_q(W_{2,a})=Y_aY_{aq^2}+Y_aY_{aq^4}^{-1}+Y_{aq^2}^{-1}Y_{aq^4}^{-1}
= Y_aY_{aq^2}\left(1 + A_{aq^3}^{-1} + A_{aq}^{-1}A_{aq^3}^{-1}\right).
\]
More generally, we have
\begin{equation}\label{eqKR}
\chi_q(W_{k,a})=\prod_{j=0}^{k-1}Y_{aq^{2j}}\left(1+A_{aq^{2k-1}}^{-1}\left(1+A_{aq^{2k-3}}^{-1}
\left(1+\cdots\left(1+A_{aq}^{-1}\right)\cdots \right)\right)\right).
\end{equation}

Following Chari and Pressley \cite{CP2}, we now describe the $q$-characters
of all the simple objects of $\CC$.
We call $q$-segment of origin $a$ and length $k$
the string of complex numbers 
\[
\Si(k,a)=\{a,\ aq^2,\ \ldots, aq^{2k-2}\}.
\]
Two $q$-segments are said to be in special position if one does 
not contain the other, and their union is a $q$-segment.
Otherwise we say that they are in general position.
It is easy to check that every finite multi-set 
$\{b_1,\ldots, b_s\}$ of elements of $\C^*$ 
can be written uniquely as a union of $q$-segments $\Si(k_i,a_i)$ 
in such a way that
every pair $(\Si(k_i,a_i),\,\Si(k_j,a_j))$ is in general position. 
Then, Chari and Pressley have proved that the simple 
module $S$
with highest {\em l}-weight 
\[
m_S =\prod_{j=1}^s Y_{b_j}
\] 
is isomorphic to the tensor product of Kirillov-Reshetikhin modules $\bigotimes_i W_{k_i,a_i}$.
Hence $\chi_q(S)$ can be calculated using (\ref{eqKR}).

\subsection{Trigonometric solutions of the Yang-Baxter equation}

Let us briefly indicate how quantum loop algebras give rise to families
of solutions of the quantum Yang-Baxter equation. A nice introduction to
these ideas is given by Jimbo in \cite{J2}. 

Consider the tensor product $W_{k,a}\otimes W_{k,b}$ of Kirillov-Reshetikhin
modules for $U_q(L\Sl_2)$. For a generic choice of $u=a/b\in\C^*$, the 
$q$-segments $\Si(k,a)$ and $\Si(k,b)$ are in general position, and therefore
$W_{k,a}\otimes W_{k,b}$ is irreducible. Moreover, since the Grothendieck
group is commutative, when the tensor product is irreducible it is isomorphic to
$W_{k,b}\otimes W_{k,a}$. Therefore, there exists up to normalization 
a unique isomorphism 
\[
I(a,b) :  W_{k,a}\otimes W_{k,b} \stackrel{\sim}{\longrightarrow} W_{k,b}\otimes W_{k,a}.
\]
Now, if $\Si(k,c)$ is another $q$-segment in general position with
$\Si(k,a)$ and $\Si(k,b)$, then  
$I_{12}(b,c)I_{23}(a,c)I_{12}(a,b)$
and 
$I_{23}(a,b)I_{12}(a,c)I_{23}(b,c)$
are two isomorphisms between the irreducible modules
$W_{k,a}\otimes W_{k,b}\otimes W_{k,c}$ and
$W_{k,c}\otimes W_{k,b}\otimes W_{k,a}$, hence they are proportional.
These intertwinners can be normalized in such a way that 
\[
I_{12}(b,c)I_{23}(a,c)I_{12}(a,b) =  I_{23}(a,b)I_{12}(a,c)I_{23}(b,c).
\]
Putting $R(a,b) = P\cdot I(a,b)$ where $P$ is the linear map from
$W_{k,b}\otimes W_{k,a}$ to $W_{k,a}\otimes W_{k,b}$ defined by
$P(w\otimes w') = w'\otimes w$, it follows that 
\[
R_{23}(b,c)R_{13}(a,c)R_{12}(a,b) =  R_{12}(a,b)R_{13}(a,c)R_{23}(b,c).
\]
Moreover it can be seen that $R(a,b)$ only depends on $a/b$, thus
setting $u=a/b$ and $v=b/c$, we obtain that
\[
R_{23}(v)R_{13}(uv)R_{12}(u) =  R_{12}(u)R_{13}(uv)R_{23}(v),
\] 
that is, $R(u)$ is a solution of the Yang-Baxter equation (\ref{YB}).
These solutions were obtained by Tarasov \cite{T1,T2}.
For example, the solution coming from the 2-dimensional representation
$W_{1,a}$ can be written in matrix form as
\[
R(u) =
\left( 
\begin{matrix}
1& 0 & 0 & 0 \\[2mm]
0& \ds\frac{q(u-1)}{u-{q}^2} & \ds\frac{1-{q}^2}{u-q^2} & 0\\[3mm]
0& \ds\frac{u(1-{q}^2)}{u-{q}^2} & \ds\frac{{q}(u-1)}{u-{q}^2} & 0 \\[3mm]
0& 0 & 0 & 1
\end{matrix}
\right). 
\]
This is the $R$-matrix associated with two famous integrable models:
the spin 1/2 XXZ chain, and the six-vertex model. 
Note that it is well-defined and invertible if and only if 
$u\not = q^{\pm 2}$. In fact, for $u = q^{\pm 2}$, Eq.~(\ref{eqKR1})
shows that the tensor product $W_{1,au}\otimes W_{1,a}$ is not irreducible, and one can
check that $W_{1,au}\otimes W_{1,a}$ is not isomorphic to
$W_{1,a}\otimes W_{1,au}$. 

More generally, the same method can be applied to any finite-dimensional
irreducible representation $W$ of $U_q(L\gg)$, using the general fact
that $W(u)\otimes W$ is irreducible except for a finite number of values
of $u\in\C^*$. We shall return to this special feature of quantum
loop algebras in \S\ref{sect_4}.

\subsection{Subcategories}\label{subsect_2.5}

Since the Dynkin diagram of $\gg$ is a bipartite graph, we have a 
partition $I=I_0\sqcup I_1$ such that every edge connects
a vertex of $I_0$ with a vertex of $I_1$.
For $i\in I$ we set
\begin{equation}
\xi_i = 
\left\{
\begin{array}{ll}
0& \mbox{if $i\in I_0$,}\\
1& \mbox{if $i\in I_1$,}
\end{array}
\right..
\end{equation}
Let $\M_\Z$ be the subset of $\M$ consisting of all monomials
in the variables
\[
 Y_{i,q^{\xi_i + 2k}},\qquad (i\in I,\ k\in\Z).
\]
Let $\CC_\Z$ be the full subcategory of $\CC$ whose objects 
$V$ have all their composition factors of the form $L(m)$
with $m\in\M_\Z$.
One can show that $\CC_\Z$ is an abelian subcategory, stable
under tensor products.
Its Grothendieck ring $R_\Z$ is the subring of $R$
generated by the classes of the fundamental modules
\[
L(Y_{i,q^{2k+\xi_i}})\qquad (i\in I,\ k\in\Z).
\]
It is known that every simple object $S$ of $\CC$ can be written 
as a tensor product $S_1(a_1)\otimes\cdots\otimes S_k(a_k)$ 
for some simple objects $S_1, \ldots, S_k$ of $\CC_\Z$ and
some complex numbers $a_1,\ldots,a_k$ such that
\[
\frac{a_i}{a_j} \not \in q^{2\Z}, \qquad (1\le i < j \le k).
\]
(Here $S_j(a_j)$ denotes the image
of $S_j$ under the auto-equivalence $\tau_{a_j}^*$, see \S\ref{subsect_2.1}.)
Therefore, the description of the simple objects of $\CC$
essentially reduces to the description of the simple
objects of $\CC_\Z$.

We will now introduce, following \cite{HL}, an increasing sequence 
of subcategories of $\CC_\Z$. 
Let $\ell\in\N$. 
Let $\M_{\ell}$ be the subset of $\M_\Z$ consisting of all monomials
in the variables
\[
 Y_{i,q^{\xi_i + 2k}},\qquad (i\in I,\ 0\le k\le \ell).
\]
Define $\CC_\ell$ to be the full subcategory of $\CC$ whose objects 
$V$ have all their composition factors of the form $L(m)$
with $m\in\M_\ell$.
\begin{proposition}[\cite{HL}]\label{propCCl}
$\CC_\ell$ is an abelian monoidal category, with Grothendieck ring
the polynomial ring
\[
 R_\ell = \Z\left[[L(Y_{i,q^{2k+\xi_i}})];\ i\in I,\ 0\le k\le \ell\right].
\]
\end{proposition}

The simple objects of the category $\CC_0$ are easy to describe.
Indeed, it follows from \cite[Prop. 6.15]{FM} that every simple
object of $\CC_0$ is a product
of fundamental modules of $\CC_0$, and conversely
any tensor product of fundamental modules of $\CC_0$ is simple.
We will see in \S\ref{sect_4} that the simple objects of the subcategory $\CC_1$ 
are already non trivial,
and that they have a nice description involving cluster
algebras.

Clearly, every simple object of $\CC_\Z$ is of the 
form $S(q^k)$ for some $k\in\Z$ and some simple object $S$ in $\CC_\ell$ 
with $\ell$ large enough.
Therefore, the description of the simple objects of $\CC$
eventually reduces to the description of the simple
objects of $\CC_\ell$ for arbitrary $\ell\in\N$.


\section{Nakajima quiver varieties and irreducible $q$-characters}
\label{sect_3}

The characters of the irreducible finite-dimensional $U_q(\gg)$-modules
are identical to those of the corresponding $\gg$-modules, and
are thus given by the classical Weyl character formula.
Moreover, Kashiwara's theory of crystal bases gave rise to a uniform
combinatorial description of these characters, generalizing the 
Young tableaux descriptions available for $\gg = \Sl_n$.
The situation is much more complicated for $U_q(L\gg)$.
Indeed, there is no analogue of Weyl's formula in this case,
and it is believed that only Kirillov-Reshetikhin modules 
(and their irreducible tensor products) have a crystal basis
(see \cite{OS} for the existence of crystals of KR-modules
for $\gg$ of classical type). 
However, inspired by earlier work on Springer theory for affine Hecke algebras, 
Ginzburg and Vasserot \cite{GV} gave a geometric description
of the irreducible $q$-characters of $U_q(L\gg)$ for $\gg=\Sl_n$ in terms
of intersection cohomology of closures of graded nilpotent orbits. 
This was extended to all simply-laced types by Nakajima \cite{N1}, using 
a graded version of his quiver varieties.
In this second lecture, we shall review Nakajima's geometric approach.

\subsection{Graded vector spaces}\label{subsect_3.1}
Recall the partition $I=I_0\sqcup I_1$ of \S\ref{subsect_2.5}.
Define the sets of ordered pairs: 
\[
\hI= I \times \Z,\quad
\hI_0=(I_0\times 2\Z) \sqcup (I_1\times (2\Z + 1)),\quad
\hI_1=(I_0\times (2\Z+1)) \sqcup (I_1\times 2\Z).
\]
We will consider finite-dimensional $\hI$-graded $\C$-vector spaces.
More precisely, we will use the letters $V, V', \ldots$ for $\hI_1$-graded
vector spaces, and the letters $W, W', \ldots $ for $\hI_0$-graded vector spaces.
We shall write
\[
V = \bigoplus_{(i,r)\in \hI_1} V_i(r),\qquad
W = \bigoplus_{(i,r)\in \hI_0} W_i(r),
\]
where the spaces $V_i(r)$, $W_i(r)$ are finite-dimensional, 
and nonzero only for a finite number of $(i,r)$.
We write $V \le V'$ if and only if $\dim V_i(r) \ge \dim V'_i(r)$
for every $(i,r)\in\hI_1$.

Consider a pair $(V,W)$ where $V$ is $\hI_1$-graded and $W$ is
$\hI_0$-graded.
We say that $(V,W)$ is $l$-dominant if
\begin{equation}\label{def_d}
d_i(r,V,W):=\dim W_i(r) - \dim V_i(r+1) - \dim V_i(r-1) - \sum_{j\not = i} a_{ij} \dim V_j(r) \ge 0 
\end{equation}
for every $(i,r)\in \hI_0$.
The pair $(0,W)$ is always $l$-dominant, and for a given $W$
there are finitely many isoclasses of $\hI_1$-graded spaces $V$ such that 
$(V,W)$ is $l$-dominant.

\subsection{ADHM equations}
Let $(V,W)$ be a pair of vector spaces, where $V$ is $\hI_1$-graded and $W$ is
$\hI_0$-graded. 
Define
\begin{align*}
&L^\bullet(V,W)= \bigoplus_{(i,r)\in \hI_1} \Hom(V_i(r),W_i(r-1)),
\\
&L^\bullet(W,V)= \bigoplus_{(i,r)\in \hI_0} \Hom(W_i(r),V_i(r-1)), 
\\
&E^\bullet(V)= \bigoplus_{(i,r)\in \hI_1;\ j:\,a_{ij}=-1} \Hom(V_i(r),V_j(r-1)).  
\end{align*}
Put  
$M^\bullet(V,W) = E^\bullet(V) \oplus L^\bullet(W,V) \oplus L^\bullet(V,W)$.
An element of $M^\bullet(V,W)$ is written 
$(B,\al,\be)$, and its components are denoted by:
\begin{align*}
&B_{ij}(r)\in \Hom(V_i(r),V_j(r-1)),
\\
&\al_i(r)\ \in \Hom(W_i(r),V_i(r-1)),
\\
&\be_i(r)\ \in \Hom(V_i(r),W_i(r-1)).
\end{align*}
We define a map 
$\mu:\ M^\bullet(V,W) \to \displaystyle\bigoplus_{(i,r)\in\hI_1} \Hom(V_i(r),V_i(r-2))$ by
\[
\mu_{(i,r)}(B,\al,\be) = \al_i(r-1)\be_i(r) + \sum_{j:\, a_{ij}=-1}  B_{ji}(r-1)B_{ij}(r),
\qquad ((i,r)\in\hI_1).
\]
We can then introduce $\L^\bullet(V,W) := \mu^{-1}(0) \subset M^\bullet(V,W)$. 
In other words, $\L^\bullet(V,W)$ is the subvariety of the affine space $M^\bullet(V,W)$
defined by the so-called complex Atiyah-Drinfeld-Hitchin-Manin equations (or ADHM, in short):
\begin{equation}\label{ADHM}
\al_i(r-1)\be_i(r) + \sum_{j:\,a_{ij}=-1}  B_{ji}(r-1)B_{ij}(r) = 0,
\qquad ((i,r)\in\hI_1).
\end{equation}

\subsection{Graded quiver varieties}
A point $(B,\al,\be)$ of $\L^\bullet(V,W)$ is called stable if the following condition
holds: for every $\hI_1$-graded subspace $V'$ of $V$, if $V'$ is $B$-invariant and contained
in $\Ker\be$ then $V'=0$.
The stable points form an open subset of $\L^\bullet(V,W)$ denoted by 
$\L^{\bullet}_s(V,W)$. 
Let 
\[
G_V := \prod_{(i,r)\in\hI_1}GL(V_i(r)).
\]
This reductive group acts on $M^\bullet(V,W)$ by base change in $V$:
\[
 g\cdot (B,\al,\be) = \left((g_j(r-1)B_{ij}(r) g_i(r)^{-1}),\ (g_i(r-1)\al_i(r)),\ (\be_i(r)g_i(r)^{-1})\right).
\]
Note that there is no action on the space $W$.
This action preserves the subvariety $\L^\bullet(V,W)$ and the open subset 
$\L^{\bullet}_s(V,W)$. Moreover, the action on $\L^{\bullet}_s(V,W)$ is
free. One can then define, following Nakajima,
\[
\fM^\bullet(V,W) := \L^{\bullet}_s(V,W) \slash G_V.
\]
This set-theoretic quotient coincides with a quotient in the geometric invariant
theory sense.
The $G_V$-orbit through $(B,\al,\be)$, considered as a point of $\fM^\bullet(V,W)$,
will be denoted by $[B,\al,\be]$. 
Note that $\fM^\bullet(V,W)$ may be empty (if there is no stable point).

One also defines the affine quotient
\[
\fM^\bullet_0(V,W) := \L^{\bullet}(V,W) \sslash G_V.
\]
By definition, the coordinate ring of $\fM^\bullet_0(V,W)$
is the ring of $G_V$-invariant functions on $\L^{\bullet}(V,W)$,
and $\fM^\bullet_0(V,W)$ parametrizes the closed $G_V$-orbits.
Since the orbit $\{0\}$ is always closed, $\fM^\bullet_0(V,W)$
is never empty. 
We have a projective morphism
\[
 \pi_V : \fM^\bullet(V,W) \to \fM^\bullet_0(V,W),
\]
mapping the orbit $[B,\al,\be]$ to the unique closed $G_V$-orbit in its
closure.
Finally, the third quiver variety is
\[ 
\fL^\bullet(V,W) := \pi_V^{-1}(0).
\]

\subsection{Properties}\label{subsect_3.4}

If it is not empty, the variety $\fM^\bullet(V,W)$ is smooth of dimension
\[
\dim \fM^\bullet(V,W) = \sum_{(i,r)\in \hI_0}
\dim V_i(r+1) d_i(r,V,W) + \dim W_i(r)\dim V_i(r-1)
\]
where $d_i(r,V,W)$ is defined by (\ref{def_d}).

The coordinate ring of $\fM^\bullet_0(V,W)$ is generated by the 
following $G_V$-invariant functions on $\L^\bullet(V,W)$:
\begin{equation}\label{gen_f}
(B,\al,\be) \mapsto \<\psi \ ,\ \be_j(r-n-1)B_{j_{n-1}j}(r-n)\cdots B_{j_1j_2}(r-2)B_{ij_1}(r-1)\al_i(r)\>,
\end{equation}
where $(i,r)\in\hI_0$, $(i,j_1,j_2,\ldots,j_{n-1},j)$ is a path (possibly of length 0)
on the (unoriented) Dynkin diagram, 
and $\psi$ is a linear form on $\Hom(W_i(r),W_j(r-n-2))$ \cite{L}.
When the path is the trivial path at vertex $i$,
the function (\ref{gen_f})~is 
\[
(B,\al,\be) \mapsto \<\psi \ ,\ \be_i(r-2)\al_i(r)\>, 
\]
where $\psi$ is a linear form on $\Hom(W_i(r),W_i(r-2))$.

In particular, if $W=W_i(r)$ is supported on a single vertex $(i,r)\in\hI_0$
then every function of the form (\ref{gen_f}) is equal to 0, so
for every $V$
the coordinate ring of $\fM^\bullet_0(V,W)$ is $\C$, and
$\fM^\bullet_0(V,W)= \{0\}$.
It follows that $\fM^\bullet(V,W) = \fL^\bullet(V,W)$ in this case.

Let $\fM^{\bullet\, {\rm reg}} _0(V,W)$ be the subset of $\fM^\bullet_0(V,W)$
consisting of the closed {\em free} $G_V$-orbits. This is open, but possibly empty.
In fact, $\fM^{\bullet\, {\rm reg}} _0(V,W) \not = \emptyset$ if and only if
$\fM^\bullet(V,W) \not = \emptyset$ and the pair $(V,W)$ is $l$-dominant.
In this case, the restriction of $\pi_V$ to 
$\pi_V^{-1}(\fM^{\bullet\, {\rm reg}} _0(V,W))$ is an isomorphism, and
in particular $\fM^{\bullet\, {\rm reg}} _0(V,W)$ is non singular 
of dimension equal to $\dim \fM^\bullet(V,W)$. 

If $V \ge V'$, that is, if $V_i(r) \subseteq V'_i(r)$ for every 
$(i,r)\in \hI_1$, then we have a natural
closed embedding $\fM^\bullet_0(V,W) \subset \fM^\bullet_0(V',W)$.
One defines
\[
\fM^\bullet_0(W) = \bigcup_V \fM^\bullet_0(V,W). 
\]
In fact, one has a stratification
\[
\fM^\bullet_0(W) = \bigsqcup_{[V]} \fM^{\bullet\,{\rm reg}}_0(V,W), 
\]
where $V$ runs through the $\hI_1$-graded spaces such that 
$(V,W)$ is $l$-dominant, and $[V]$ denotes the isomorphism class
of $V$ as a graded space.
It follows from \S\ref{subsect_3.1} that 
$\fM^\bullet_0(W)$ has finitely many strata.

\subsection{Examples}\label{subsect_3.5}
\textbf{1.}\ 
Take $\gg=\Sl_2$ of type $A_1$.
Assume that $I_0 = I = \{1\}$. Since $I$ is a singleton, we can 
drop indices $i$ in the notation and write $\hI_0 = 2\Z$,
$\hI_1 = 2\Z+1$.
Hence 
\[
W = \bigoplus_{r\in 2\Z} W(r), 
\qquad
V = \bigoplus_{s\in 2\Z+1} V(s), 
\]
and $M^\bullet(V,W) = L^\bullet(W,V) \oplus L^\bullet(V,W)$
consists of pairs $(\al,\be)$: the $B$-component is zero
in this case. In particular, any subspace $V'$ of $V$ is $B$-stable, so
$(\al,\be)$ is stable if and only $\be(s)$ is injective
for every $s\in 2\Z+1$. The ADHM equations reduce to
\[
\al(s-1)\be(s) = 0,
\qquad (s\in 2\Z+1).
\]
With any pair $(\al,\be)$ we associate $x = \be\al \in \End(W)$, 
and $E=\be(V) \subseteq W$.
Clearly, $x$ and $E$ depend only on the $G_V$-orbit of $(\al,\be)$.

Let $\NN(W)$ denote the subvariety of $\End(W)$ consisting of 
degree -2 endomorphisms of $W$ (that is,
$x(W(r)) \subseteq W(r-2)$ for every $r$) satisfying $x^2 = 0$.
Let $\F(V,W)$ denote the  
variety of pairs $(x,E)$, where $E$ is a graded subspace of $W$
with $\dim E(r) = \dim V(r+1)$, and $x\in \NN(W)$ is such that
\[
 \im x \subseteq E \subseteq  \Ker x.
\]
Then one can check that the map $[\al,\be]\mapsto(x,E)$
establishes an isomorphism 
\[
\fM^\bullet(V,W) \stackrel{\sim}\longrightarrow \F(V,W).
\]
Moreover, the affine variety $\fM^\bullet_0(W)$ is isomorphic to $\NN(W)$,
and the projective morphism $\pi_V : \fM^\bullet(V,W) \to \fM^\bullet_0(V,W) \subseteq \fM^\bullet_0(W)$
is the first projection 
\[
\pi_V(x,E) = x.
\] 
Thus the fiber $\fL^\bullet(V,W) = \pi_V^{-1}(0)$ is the Grassmannian 
of graded subspaces $E$ of $W$ with $\dim E(r) = \dim V(r+1)$.
Finally, $\fM^\bullet_0(V,W)$ is isomorphic to the subvariety of $\NN(W)$
defined by the rank conditions:
\[
\dim x(W(r)) \le \dim V(r-1), \qquad (r\in 2\Z). 
\]

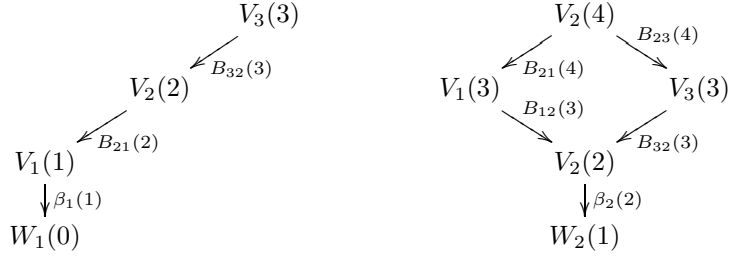
\begin{figure}[t]
\[
\xymatrix@-1.0pc{
&&
&\ar[ld]^{B_{32}(3)} V_3(3) 
\\
&&\ar[ld]^{B_{21}(2)} V_2(2) &&
\\
&V_1(1)\ar[d]^{\be_1(1)}  && 
\\
&W_1(0)
}
\quad
\xymatrix@-1.0pc{
&&\ar[ld]^{B_{21}(4)} V_2(4) \ar[rd]^{B_{23}(4)}&&
\\
&{V_1(3)}\ar[rd]^{B_{12}(3)}&
&\ar[ld]^{B_{32}(3)} V_3(3) 
\\
&& V_2(2)\ar[d]^{\be_2(2)} &&
\\
&&W_2(1)
}
\]
\caption{\label{Fig1} {\it The cases $W=W_1(0)=\C$, and $W=W_2(1)=\C$, in type $A_3$.}}
\end{figure}

\noindent
\textbf{2.}\ 
Take $Q$ of type $A_3$, with sink set $I_0=\{1,3\}$ and 
source set $I_1=\{2\}$.
Choose first
\[
W = W_1(0) = \C.
\] 
One can check that the only choices of $V$ such that
$\fM^\bullet(V,W) \not = \emptyset$ are 
\begin{itemize}
 \item[(a)] $V= 0$.
 \item[(b)] $V= V_1(1) = \C$.
 \item[(c)] $V= V_1(1) \oplus V_2(2)= \C \oplus \C$.
 \item[(d)] $V= V_1(1) \oplus V_2(2) \oplus V_3(3) = \C \oplus \C \oplus \C$.
\end{itemize}
For instance, this can be checked by using the dimension formula
of \S\ref{subsect_3.4} : in all other cases, this formula produces a 
negative number.
Case (d) is illustrated in Figure~\ref{Fig1}.
Let us determine $\L^\bullet_s(V,W)$ in each case.
In cases (b), (c), (d), the map $\be_1(1): V_1(1) \to W_1(0)$ has to be
injective. Indeed its kernel is $B$-invariant, and so the stability condition
forces it to be trivial.
Similarly, in cases (c), (d), the map $B_{21}(2): V_2(2) \to V_1(1)$ has to be
injective. Indeed its kernel is $B$-invariant and contained in $\Ker\be_2(2) = V_2(2)$, 
and so the stability condition forces it again to be trivial.
Finally, in case (d), for the same reasons, the map $B_{32}(3): V_3(3) \to V_2(2)$ 
has to be injective.
Thus we have
\begin{itemize}
 \item[(b)] $\L^\bullet_s(V,W)=\C^*$ and $\fM^\bullet(V,W)=\C^*/\C^* = \{{\rm point}\}$.
 \item[(c)] $\L^\bullet_s(V,W)=\C^*\times \C^*$ and 
$\fM^\bullet(V,W)=(\C^*\times \C^*)/(\C^*\times \C^*) = \{{\rm pt}\}$.
 \item[(d)] $\L^\bullet_s(V,W)=\C^*\times \C^*\times \C^*$,  
$\fM^\bullet(V,W)=(\C^*\times \C^*\times \C^*)/(\C^*\times \C^*\times \C^*) = \{{\rm pt}\}$.
\end{itemize} 
As a second choice, take
\[
W = W_2(1) = \C.
\] 
One can check that the only choices of $V$ such that
$\fM^\bullet(V,W) \not = \emptyset$ are 
\begin{itemize}
 \item[(a)] $V= 0$.
 \item[(b)] $V= V_2(2) = \C$.
 \item[(c)] $V= V_2(2) \oplus V_1(3)= \C \oplus \C$.
 \item[(d)] $V= V_2(2) \oplus V_3(3)= \C \oplus \C$.
 \item[(e)] $V= V_2(2) \oplus V_1(3) \oplus V_3(3) = \C \oplus \C \oplus \C$.
 \item[(f)] $V= V_2(2) \oplus V_1(3) \oplus V_3(3) \oplus V_2(4) = \C \oplus \C \oplus \C \oplus \C$.
\end{itemize}
Case (f) is illustrated in Figure~\ref{Fig1}.
Let us determine $\L^\bullet_s(V,W)$ in each case.
In cases (b), (c), (d), (e) the map $\be_2(2): V_2(2) \to W_2(1)$ has to be
injective. Indeed its kernel is $B$-invariant, and so the stability condition
forces it to be trivial.
Similarly, in cases (c), (e), (f) the map $B_{12}(3)$ has to be
injective. Indeed its kernel is $B$-invariant and contained in $\Ker\be_1(3) = V_1(3)$, 
and so the stability condition forces it again to be trivial.
Similarly, in cases (d), (e), (f) the map $B_{32}(3)$ has to be
injective.
Finally, in case (f), the ADHM equations imply the relation
\begin{equation}\label{rel}
B_{12}(3)B_{21}(4) + B_{32}(3)B_{23}(4) = 0.
\end{equation}
Since $B_{12}(3)$ and $B_{32}(3)$ are injective, this implies that
$B_{21}(4)$ and $B_{23}(4)$ are both injective or both equal to 0.
If they are both equal to 0 then $V_2(4)$ is $B$-invariant and
contained in $\Ker\be_2(4)$, so $(B,\al,\be)$ is not stable. Hence
$B_{21}(4)$ and $B_{23}(4)$ are both injective.
Thus we have
\begin{itemize}
 \item[(b)] $\L^\bullet_s(V,W)=\C^*$ and $\fM^\bullet(V,W)=\C^*/\C^* = \{{\rm pt}\}$.
 \item[(c)] $\L^\bullet_s(V,W)=\C^*\times \C^*$ and 
$\fM^\bullet(V,W)=(\C^*\times \C^*)/(\C^*\times \C^*) = \{{\rm pt}\}$.
\item[(d)] $\L^\bullet_s(V,W)=\C^*\times \C^*$ and 
$\fM^\bullet(V,W)=(\C^*\times \C^*)/(\C^*\times \C^*) = \{{\rm pt}\}$.
 \item[(e)] $\L^\bullet_s(V,W)=\C^*\times \C^*\times \C^*$, 
$\fM^\bullet(V,W)=(\C^*\times \C^*\times \C^*)/(\C^*\times \C^*\times \C^*) = \{{\rm pt}\}$.
 \item[(f)] $\L^\bullet_s(V,W)=\C^*\times \C^*\times \C^*\times \C^*$, because $B_{23}(4)$
can be expressed in terms of $B_{12}(3), B_{21}(4)$, and
$B_{32}(3)$ in view of (\ref{rel}). Hence, we find again that
$\fM^\bullet(V,W)=(\C^*\times \C^*\times \C^* \times \C^*)/(\C^*\times \C^*\times \C^*\times \C^*) = \{{\rm pt}\}$.
\end{itemize}

\noindent
\textbf{3.}\
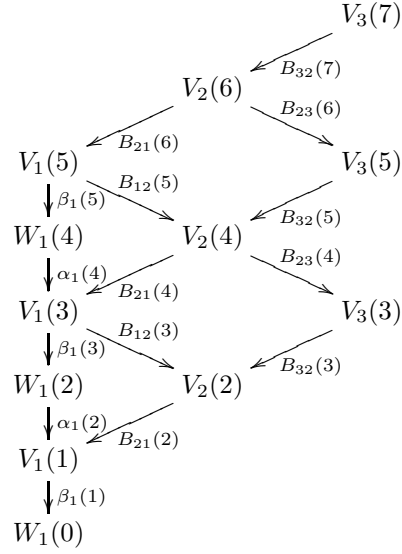
\begin{figure}[t]
\[
\xymatrix@-1.0pc{
&\qquad\qquad\qquad&&
\\
  &&  \ar[ld]^{B_{32}(7)} V_3(7)
\\
&\ar[ld]^{B_{21}(6)} V_2(6)\ar[rd]^{B_{23}(6)} 
\\
\ar[d]^{\be_1(5)}V_1(5)\ar[rd]^{B_{12}(5)}&&\ar[ld]^{B_{32}(5)} V_3(5)
\\
\ar[d]^{\al_1(4)}W_1(4)&\ar[ld]^{B_{21}(4)} V_2(4)\ar[rd]^{B_{23}(4)} 
\\
\ar[d]^{\be_1(3)}V_1(3)\ar[rd]^{B_{12}(3)}&&\ar[ld]^{B_{32}(3)} V_3(3)
\\
\ar[d]^{\al_1(2)} W_1(2)&\ar[ld]^{B_{21}(2)} V_2(2) &&
\\
V_1(1)\ar[d]^{\be_1(1)}  && 
\\
W_1(0)
}
\]
\caption{\label{Fig2} {\it The case $W=W_1(0)\oplus W_1(2)\oplus W_1(4)$ in type $A_3$.}}
\end{figure}
Assume that $Q$ is of type $A_n$ and $1\in I_0$.
Take $W$ of the form 
\[
W = W_1(r) \oplus W_1(r+2) \oplus \cdots \oplus W_1(r+2k),
\qquad
(r\in 2\Z,\ k\in\N),
\]
where $\dim W_1(r+2j)=d^{(j)}\ (0\le j\le k)$.
Thus, $W$ can be regarded as a $2\Z$-graded vector space.
One can check that the map 
\[
(B,\al,\be)\mapsto x:=(\be_1(r+1)\al_1(r),\be_1(r+3)\al_1(r+2)\ldots,\be_1(r+2k-1)\al_1(r+2k))
\]
induces an isomorphism from $\fM^\bullet_0(W)$ to the variety 
of degree -2 endomorphisms $x$ of $W$
satisfying $x^{n+1}=0$.
In other words, $\fM^\bullet_0(W)$ is the affine space of representations
in $W$ of the quiver
\[
r \stackrel{\ga_1}{\leftarrow} r+2 \stackrel{\ga_2}{\leftarrow} r+4 
\stackrel{\ga_3}{\leftarrow} \cdots \stackrel{\ga_k}{\leftarrow} r+2k 
\]
bound by the relations 
\[
\ga_i\ga_{i+1}\cdots\ga_{i+n} = 0,\qquad (1\le i \le k-n).
\] 
Let us now determine $\L^\bullet_s(V,W)$ for a given $\hI_1$-graded space $V$.
First, the stability condition implies that 
$V_i(j) = 0$ for $j < r+i$.
Next, it is easy to show by induction that if $(B,\al,\be)$ is stable then
the following maps are injective:
\[
 \be_1(r+1+2j)\quad (0\le j\le k),\qquad B_{i,i-1}(r+i+2j)\quad (2\le i\le n,\ 0\le j\le k).
\]
Moreover, we have $V_i(j) = 0$ for $j > r+i+2k$.
Therefore a typical example of $(B,\al,\be)$ in $\L_s^\bullet(V,W)$ looks like in Figure~\ref{Fig2} with
all maps $\be_1(j)$ and $B_{i,i-1}(j)$ injective.
Put 
\[
E_1^{(j)} := \be_1(r+1+2j)(V_1(r+1+2j)) \subseteq W_1(r+2j), \qquad (0\le j\le k),
\]
and for $i=2,\ldots,n$,
\[
E_i^{(j)} := \be_1(r+1+2j)B_{21}(r+2+2j)\cdots B_{i,i-1}(r+i+2j)(V_i(r+i+2j)).
\]
The vector spaces $E_i=\bigoplus_{j} E_i^{(j)}$ form an $n$-step flag 
\[
F^\bullet = \left(E_0 = W\supseteq E_1 \supseteq \cdots \supseteq E_n \supseteq 0 = E_{n+1}\right) 
\]
of graded subspaces of $W = E_0$, with graded dimension
\[
\bd_i = (\dim V_i(r+i),\, \dim V_i(r+i+2),\, \ldots),
\qquad (i=1,\ldots n). 
\]
We get a well-defined map $(B,\al,\be) \mapsto (x,F^\bullet)$ 
from $\L_s^\bullet(V,W)$ to the variety $\F(V,W)$ of pairs consisting of a graded
flag $F^\bullet$ of dimension $(\bd_1,\ldots,\bd_n)$ in $W$, together with a graded
nilpotent endomorphism $x$ preserving this flag (that is, such that 
$x(E_i^{(j)})\subseteq E_{i+1}^{(j-1)}$).
This map is $G_V$-equivariant, hence induces a map 
$\fM^\bullet(V,W) \to \F(V,W)$, and one can check that this
is an isomorphism. Moreover, in this identification the map 
$\pi_V : \fM^\bullet(V,W) \to \fM_0^\bullet(V,W)$ becomes
the projection $(x,F^\bullet) \mapsto x$.
Finally, the zero fibre $\fL^\bullet(V,W)=\pi_V^{-1}(0)$ is isomorphic
to the variety of
graded flags of $W$ of dimension $(\bd_1,\ldots,\bd_n)$.

More generally, the fiber $\pi_V^{-1}(x)$ is the variety of 
graded flags of $W$ of dimension $(\bd_1,\ldots,\bd_n)$
which are preserved by $x$. Ginzburg and Vasserot \cite{GV}
have shown that the Borel-Moore homologies of the varieties
\[ 
M_x = \bigsqcup_{V} \pi_V^{-1}(x),\qquad (x\in \fM_0^\bullet(W)),
\]
(where $V$ runs over isoclasses of $\hI_1$-graded spaces)
have natural structures of $U_q(\hSl_{n+1})$-modules, called
standard modules.
These modules are not simple, but can be decomposed into simple
ones using the decomposition theorem for perverse sheaves. 
In the next sections, we review Nakajima's extension of these
results to other root systems.

\subsection{Quiver varieties and standard modules for $U_q(L\gg)$}

To a pair $(V,W)$ of graded spaces as in \S\ref{subsect_3.1} we attach
two monomials in $\Y$ given by
\begin{equation}
Y^W = \prod_{(i,\,r)\in \hI_0} Y_{i,q^r}^{\dim W_i(r)},
\qquad
A^V = \prod_{(j,\,s)\in \hI_1} A_{j,q^s}^{-\dim V_j(s)}.
\end{equation}
One can check that the monomial $Y^WA^V$ is dominant if and only
if the pair $(V,W)$ is {\em l}-dominant in the sense of \S\ref{subsect_3.1}.

Let us associate with $W$ the simple $U_q(L\gg)$-module
$L(W):=L(Y^W)$, which belongs to the subcategory $\CC_\Z$.
We can also attach to $W$ the tensor product 
\[
M(W) = \bigotimes_{(i,\,r)\in \hI_0} L(Y_{i,q^r})^{\otimes \dim W_i(r)}.   
\]
This product is not simple in general. Moreover, its isomorphism class
may depend on the chosen ordering of the factors. However, we will only
be interested in its $q$-character (or in its class in $R_\Z$) which
is independent of this ordering. The modules $M(W)$ are called standard
modules.

The morphism $\pi_V : \fM^\bullet(V,W) \to \fM^\bullet_0(V,W)$ being projective,
its zero fiber $\fL^\bullet(V,W)$ is a complex
projective variety. Let $\chi(\fL^\bullet(V,W))$ denote
its Euler characteristic. Note that $\fL^\bullet(V,W)$ 
has no odd cohomology \cite[\S7]{N1}, so $\chi(\fL^\bullet(V,W))$
is equal to the total dimension of the cohomology.
\begin{theorem}[Nakajima \cite{N1}]\label{theo_std}
The $q$-character of the standard module $M(W)$ is given by
\[
\chi_q(M(W)) = Y^W \sum_{[V]} \chi(\fL^\bullet(V,W))\, A^V, 
\]
where the sum runs over all isomorphism classes $[V]$ of
$\hI_1$-graded spaces $V$. 
\end{theorem}
If $W$ has dimension 1, $M(W)=L(W)$ is a fundamental
module, hence Theorem~\ref{theo_std} describes in particular the
$q$-characters of all fundamental modules. On the other hand,
$\chi_q(M(W))$ is the product of the $q$-characters of the factors
of $M(W)$, so it can also be expressed as a product of $q$-characters of fundamental modules.

\subsection{Examples}\label{subsect_3.7}
We illustrate Theorem~\ref{theo_std} with the examples of 
\S\ref{subsect_3.5}.

\medskip
\noindent
\textbf{1.}\ 
We take $\gg = \Sl_2$. 
The variety $\fL^\bullet(V,W)$ is isomorphic to the 
product of ordinary Grassmannians
\[
\prod_{r\in 2\Z} \Gr\left(\dim V(r+1), \dim W(r)\right), 
\]
so its Euler characteristic is
\[
\prod_{r\in 2\Z} 
\left(
\begin{matrix}
\dim W(r) \\[1.5mm]
\dim V(r+1) 
\end{matrix}
\right).
\]
Hence 
\[
\chi_q(M(W)) = Y^W \sum_{[V]}  
\prod_{r\in 2\Z}
\left(
\begin{matrix}
\dim W(r) \\[1.5mm]
\dim V(r+1) 
\end{matrix}
\right)
A^V.
\]
On the other hand, recall that $\chi_q(L(Y_{q^r})) = Y_{q^r}(1 + A_{q^{r+1}}^{-1})$. 
We can check that the above value of $\chi_q(M(W))$ is equal to
\[
 \prod_{r\in 2\Z} \chi_q(L(Y_{q^r}))^{\dim W(r)},
\]
as it should.

\medskip
\noindent
\textbf{2.}\ 
Take $\gg$ of type $A_3$.
It follows from the calculations of \S\ref{subsect_3.5}{.2}
that 
\begin{align*}
\chi_q(L(Y_{1,q^0}))&= Y_{1,q^0}\left(1 + A_{1,q}^{-1} + A_{1,q}^{-1}A_{2,q^2}^{-1}  
+ A_{1,q}^{-1}A_{2,q^2}^{-1} A_{3,q^3}^{-1}\right),
\\
\chi_q(L(Y_{2,q^1}))&= Y_{2,q^1}\left(1 + A_{2,q^2}^{-1} + A_{2,q^2}^{-1}A_{1,q^3}^{-1} 
+ A_{2,q^2}^{-1}A_{3,q^3}^{-1} + A_{2,q^2}^{-1}A_{1,q^3}^{-1}A_{3,q^3}^{-1} \right.
\\ 
&\ \ \ \left.+\ A_{2,q^2}^{-1}A_{1,q^3}^{-1} A_{3,q^3}^{-1}A_{2,q^4}^{-1}\right).
\end{align*}

\medskip
\noindent
\textbf{3.}\ 
Assume that $\gg$ is of type $A_n$.
Choosing $W$ of dimension 1, it follows from \S\ref{subsect_3.5}{.3} that
\begin{align*}
\chi_q(L(Y_{1,q^r}))&= Y_{1,q^r}\left(1 + A_{1,q^{r+1}}^{-1} + A_{1,q^{r+1}}^{-1}A_{2,q^{r+2}}^{-1}  
\right. 
\\
&\ \ \ \left.+\ \cdots\ +\ A_{1,q^r}^{-1}A_{2,q^{r+2}}^{-1}\cdots A_{n,q^{r+n}}^{-1}\right).
\end{align*}
 
\subsection{Standard modules and the graded preprojective algebra}
\begin{figure}[t]
\[
\def\objectstyle{\scriptstyle}
\def\lablestyle{\scriptstyle}
\xymatrix@-1.0pc{
&&&&\\
&{(1,4)}\ar[rd]&{}\save[]+<0cm,2ex>*{\vdots}\restore
&\ar[ld] (3,4) 
\\
&&\ar[ld] (2,3) \ar[rd]&&
\\
&{(1,2)}\ar[rd]&
&\ar[ld] (3,2) 
\\
&&\ar[ld] (2,1) \ar[rd]&&
\\
&(1,0) \ar[rd] &&\ar[ld] (3,0)
\\
&&\ar[ld] (2,-1) \ar[rd]&&
\\
&(1,-2) &{}\save[]+<0cm,-2ex>*{\vdots}\restore& (3,-2) 
\\
}
\]
\caption{\label{Fig3} {\it The quiver $\Z Q$ in type $A_3$.}}
\end{figure}
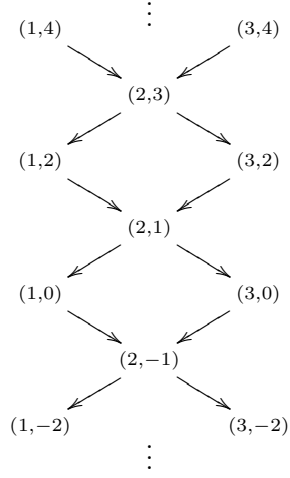
Let $Q$ be the sink-source orientation of the Dynkin diagram
of $\gg$, with set of sinks $I_0$ and set of sources $I_1$.
Define the repetition quiver $\Z Q$ as the infinite quiver 
with set of vertices $\hI_0=(I_0\times 2\Z) \sqcup (I_1\times (2\Z + 1))$,
and two types of arrows:
\begin{itemize}
 \item[(i)] for every arrow $i\to j$ in $Q$ we have 
arrows $(i,2m+1)\to (j,2m)$ in $\Z Q$ for all $m\in\Z$;
 \item[(ii)] for every arrow $i\to j$ in $Q$ we have 
arrows $(j,2m)\to (i,2m-1)$ in $\Z Q$ for all $m\in\Z$.
\end{itemize}
As an example, the quiver $\Z Q$ for $Q$ of type $A_3$ is shown
in Figure~\ref{Fig3}. 

We then introduce a set of degree two elements in the path
algebra of $\Z Q$: for every $(i,r)\in\hI_0$, 
let $\sigma_{i,r}$ be the sum of all paths from $(i,r)$
to $(i,r-2)$. The graded preprojective algebra of $Q$ is
by definition the quotient $\hL$ of the path algebra of $\Z Q$
by the two-sided ideal generated by the $\sigma_{i,r}\ ((i,r) \in \hI_0)$.
The algebra $\hL$ is well known to be the universal 
cover of the preprojective algebra $\Lambda$ of $Q$, in the
sense of \cite{G}.

It turns out that the quiver variety $\fL^\bullet(V,W)$ is
homeomorphic to a quiver Grassmannian of an injective $\hL$-module.
To state this precisely, let us denote by $S_{i,r}$ the one-dimensional 
simple $\hL$-module supported on vertex $(i,r)$ of $\Z Q$. 
Let $\De_{i,r}$ be the injective hull of $S_{i,r}$.
(This is a finite-dimensional module.) 
To the $\hI_0$-graded vector space $W$ we attach the 
injective module
\[
 \De_W = \bigoplus_{(i,r)\in \hI_0} \De_{i,r}^{\oplus \dim W_i(r)}. 
\]
To the $\hI_1$-graded vector space $V$ we attach the dimension
vector
\[
d_V = (\dim V_i(r+1);\ (i,r)\in\hI_0).
\]
We then have the following result, due to Lusztig \cite{L} in the ungraded case,
and extended to the graded case by Savage and Tingley \cite{ST}.
\begin{proposition}
The complex variety $\fL^\bullet(V,W)$ is homeomorphic to the
Grassmannian $\Gr(d_V,\De_W)$ of $\hL$-submodules of $\De_W$ with dimension vector $d_V$. 
\end{proposition}
It follows that we can rewrite Nakajima's formula for standard modules of $\CC_\Z$ as
\begin{equation}\label{eqstandard}
\chi_q(M(W)) = Y^W \sum_{[V]} \chi(\Gr(d_V,\De_W))\, A^V. 
\end{equation}
In particular, for the fundamental modules of $\CC_\Z$ we get
\begin{equation}\label{eqfund}
\chi_q(L(Y_{i,q^r})) = Y_{i,q^r} \sum_{[V]} \chi(\Gr(d_V,\De_{i,r}))\, A^V.  
\end{equation}

\begin{figure}[t]
\[
\def\objectstyle{\scriptstyle}
\def\lablestyle{\scriptstyle}
\xymatrix@-1.0pc{
&&
&\ar[lld]\ar[ld] {(3,4)} \ar[rd]
\\
&(1,3)\ar[rrd]&(2,3) \ar[rd]&&(4,3)\ar[ld]
\\
& &&\ar[lld]\ar[ld] (3,2)\ar[rd]
\\
&(1,1)\ar[rrd]& (2,1) \ar[rd]&&(4,1)\ar[ld]
\\
& && (3,0) 
\\
}
\]
\caption{\label{FigD} {\it The skeleton of the injective module $\De_{3,0}$ in type $D_4$}.}
\end{figure}
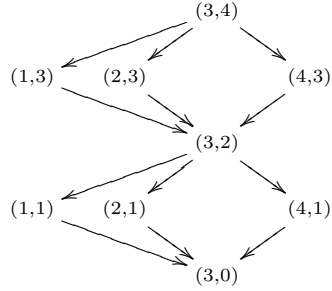

\subsection{Examples}\label{subsect_newexamples}
%
\textbf{1.}\ 
One can easily recover the formulas of \S\ref{subsect_3.7}
using the well-known description of the indecomposable injective $\hL$-modules
for a Dynkin quiver $Q$ of type $A_n$.
Thus, \S\ref{subsect_3.7}.3 follows immediately from the fact that
the injective module $\De_{1,r}$ is an $n$-dimensional
module with a unique composition series. 

More generally, all the Grassmannians of submodules of all indecomposable injective 
$\hL$-modules in type $A_n$ are reduced to a point. This implies that all 
fundamental $U_q(L\Sl_{n+1})$-modules have multiplicity-free $q$-characters,
that is, all their {\em l}-weight spaces have dimension 1.

\medskip
\noindent
\textbf{2.}\ 
Let $\gg= \so_8$ be of type $D_4$. We label by $3$ the central node of the 
Dynkin diagram, and we set $I_0 = \{3\}$, $I_1 = \{1,2,4\}$.
The injective module $\De_{3,0}$ has total dimension 10, and its dimension
vector is supported on the finite strip of $\Z Q$ displayed in
Figure~\ref{FigD}. More precisely, every vertex of the picture carries
a 1-dimensional vector space, except $(3,2)$ which has a 2-dimensional
vector space. There are 28 non-trivial quiver Grassmannians
$\Gr(d,\De_{3,0})$, and it easy to check that all of them are points,
except for the following $d$:
\[
 d_{3,0}=d_{1,1}=d_{2,1}=d_{4,1}=d_{3,2}=1,\quad
 d_{1,3}=d_{2,3}=d_{4,3}=d_{3,4}=0,
\]
for which $\Gr(d,\De_{3,0}) \simeq \P^1(\C)$. 
It follows that the fundamental module $L(Y_{3,q^0})$ has dimension 29.
More precisely, writing for short $v_{i,s}:=A_{i,q^s}^{-1}$, we have 
\begin{align*}
\chi_q(L(Y_{3,0}))=&\ Y_{3,q^0}(1 + v_{3,1} + v_{3,1}v_{1,2} + v_{3,1}v_{2,2} + v_{3,1}v_{4,2} 
+ v_{3,1}v_{1,2}v_{2,2} 
\\
&+ v_{3,1}v_{1,2}v_{4,2} + v_{3,1}v_{2,2}v_{4,2}+v_{3,1}v_{1,2}v_{2,2}v_{4,2}
+ v_{3,1}v_{1,2}v_{2,2}v_{3,3} 
\\&+ v_{3,1}v_{1,2}v_{4,2}v_{3,3} + v_{3,1}v_{2,2}v_{4,2}v_{3,3}
+ 2\, v_{3,1}v_{1,2}v_{2,2}v_{4,2}v_{3,3}
\\
&+ v_{3,1}v_{1,2}v_{2,2}v_{3,3}v_{4,4} + v_{3,1}v_{1,2}v_{4,2}v_{3,3}v_{2,4} + v_{3,1}v_{2,2}v_{4,2}v_{3,3}v_{1,4}
\\
&+v_{3,1}v_{1,2}v_{2,2}v_{4,2}v_{3,3}^2 
+ v_{3,1}v_{1,2}v_{2,2}v_{4,2}v_{3,3}v_{1,4} 
\\
&+ v_{3,1}v_{1,2}v_{2,2}v_{4,2}v_{3,3}v_{2,4} 
+v_{3,1}v_{1,2}v_{2,2}v_{4,2}v_{3,3}v_{4,4} 
\\
&+v_{3,1}v_{1,2}v_{2,2}v_{4,2}v_{3,3}^2v_{1,4}
+ v_{3,1}v_{1,2}v_{2,2}v_{4,2}v_{3,3}^2v_{2,4} 
\\
&+v_{3,1}v_{1,2}v_{2,2}v_{4,2}v_{3,3}^2v_{4,4}
+v_{3,1}v_{1,2}v_{2,2}v_{4,2}v_{3,3}^2v_{1,4}v_{2,4}
\\
&+v_{3,1}v_{1,2}v_{2,2}v_{4,2}v_{3,3}^2v_{1,4}v_{4,4}
+v_{3,1}v_{1,2}v_{2,2}v_{4,2}v_{3,3}^2v_{2,4}v_{4,4}
\\
&+v_{3,1}v_{1,2}v_{2,2}v_{4,2}v_{3,3}^2v_{1,4}v_{2,4}v_{4,4}
+v_{3,1}v_{1,2}v_{2,2}v_{4,2}v_{3,3}^2v_{1,4}v_{2,4}v_{4,4}v_{3,5})
\end{align*}
The restriction of $L(Y_{3,0})$ to $U_q(\so_8)$ decomposes into
the direct sum of the fundamental module with fundamental weight $\vpi_3$,
of dimension 28, and of a copy of the trivial representation. 

\subsection{Perverse sheaves}\label{perverse}
Recall from \S\ref{subsect_3.4} the stratification 
\[
\fM^\bullet_0(W) = \bigsqcup_{[V']} \fM^{\bullet\,{\rm reg}}_0(V',W), 
\]
where $V'$ runs through the $\hI_1$-graded spaces such that 
the pair $(V',W)$ is $l$-dominant.
We denote by $IC_W(V')$ the intersection cohomology complex associated
with the trivial local system on the stratum $\fM^{\bullet\,{\rm reg}}_0(V',W)$.

Consider an arbitrary $\hI_1$-graded space $V$ such that 
$\fM^\bullet(V,W) \not = \emptyset$.
The map 
\[
\pi_V : \fM^\bullet(V,W) \to \fM^\bullet_0(V,W)
\] 
is projective, and the variety $\fM^\bullet(V,W)$ is smooth. 
Hence, by the decomposition
theorem, the push-down $(\pi_V)_!(1_{\fM^\bullet(V,W)})$ 
of the constant sheaf on $\fM^\bullet(V,W)$ 
is a direct sum
of shifts of simple perverse sheaves in the derived category $\D(\fM^\bullet_0(V,W))$.
We can regard these perverse sheaves as objects of $\D(\fM^\bullet_0(W))$
by extending them by $0$ on the complement $\fM^\bullet_0(W) \setminus \fM^\bullet_0(V,W)$. 
Nakajima has shown that all these perverse sheaves are of the form
$IC_W(V')$ for some $l$-dominant pair $(V',W)$ with $V' \ge V$.
So we can write in the Grothendieck group of $\D(\fM^\bullet_0(W))$
\begin{equation}\label{eqIC}
[(\pi_V)_!(1_{\fM^\bullet(V,W)}[\dim\fM^\bullet(V,W)])] = \sum_{V'\ge V} a_{V,V';W}(t) [IC_W(V')]. 
\end{equation}
Here $t$ is a formal variable implementing the action of the
shift functor: 
\[
t^j[L] = [L[j]], 
\]
and $a_{V,V';W}(t) \in \N[t^{\pm 1}]$ is the graded multiplicity.
(The additional shift in degree by $\dim\fM^\bullet(V,W)$ makes the left-hand
side invariant under Verdier duality.)
Note that in (\ref{eqIC}), the pair $(V,W)$ is not necessarily
$l$-dominant.

We can now state the main result of this lecture.
\begin{theorem}[Nakajima \cite{N1}]\label{theo_sim}
Let $W$ be an $\hI_0$-graded space, and let $L(W)$ be the corresponding
simple module in $\CC_\Z$.
The coefficient of the monomial $Y^WA^V$ in $\chi_q(L(W))$ 
is equal to
$a_{V,0;W}(1)$. 
\end{theorem}
In other words, the $l$-weight multiplicities of $L(W)$ are calculated
by the (ungraded) multiplicities of the skyscraper sheaf $IC_W(0) = 1_{\{0\}}$
in the expansions of the push-downs $[(\pi_V)_!(1_{\fM^\bullet(V,W)})]$  
on the basis $\{[IC_W(V')]\}$.

\subsection{Examples}
\textbf{1.}\ 
Let $L(W) = L(Y_{i,q^r})$ be a fundamental module.
Then, by \S\ref{subsect_3.4}, $\fM^\bullet_0(W) = \{0\}$,
and $\fM^\bullet(V,W) = \fL^\bullet(V,W)$. 
So $a_{V,0;W}(t)$ is the Poincar\'e polynomial of the cohomology
of $\fL^\bullet(V,W)$ (up to some shift). Since the odd cohomology groups
of $\fL^\bullet(V,W)$ vanish, we recover that the 
coefficient of $Y^WA^V$ in $\chi_q(L(W))$ is the Euler characteristic
of $\fL^\bullet(V,W)$, in agreement with Theorem~\ref{theo_std}.

\medskip
\noindent
\textbf{2.}\ 
Take $\gg$ of type $A_1$, and
\[
W=W(r)\oplus W(r+2) \oplus \cdots \oplus W(r+2k),
\]
with $\dim W(r+2i) = 1$ for every $i=0,1,\ldots,k$.
The corresponding simple $U_q(\hg)$-module is
the Kirillov-Reshetikhin module 
$L(W) = W_{k+1,q^r}$.
Recall from \S\ref{subsect_3.5}.1 the description
of the quiver varieties $\fM^\bullet(V,W)$
and $\fM^\bullet(W)$.
The variety $\fM^\bullet(V,W)$ is non-empty if
and only if $\dim V(r+2j+1) \le \dim W(r+2j)$ for every $j=0,\ldots,k$.
For any such choice, since $\dim V(r+2j+1)$ is equal to $0$ or $1$,
there is a unique graded subspace $E$ of $W$
satisfying $\dim E(r+2j) = \dim V(r+2j+1)$.
Moreover, the set of $x$'s such that
$\im x \subseteq E \subseteq \Ker x$ is isomorphic to 
the vector space of linear maps of degree -2 from $W/E$ to $E$. 
We have two cases:

(i)\ \
If $V$ is such that there exists $s\in\{0,\ldots,k\}$ with
\[
\dim V(r+2j+1) = 
\left\{
\begin{array}{lc}
0&\mbox{if } j< s,\\
1&\mbox{if } j\ge s, 
\end{array}
\right.
\]
then this space is reduced to $\{0\}$, and 
\[
(\pi_V)_!(1_{\fM^\bullet(V,W)}) = 1_{\{0\}},
\qquad
a_{V,0,W}(1)=1.
\]

(ii)\ \ 
Otherwise, if there is $s\in\{0,\ldots,k-1\}$ with
\[
\dim V(r+2s+1)= 1,\quad \dim V(r+2s+3)=0,
\]
then this space has positive dimension, and 
$(\pi_V)_!(1_{\fM^\bullet(V,W)})$ is a simple
perverse sheaf $\not = 1_{\{0\}}$.
So $a_{V,0,W}=0$.

In conclusion, the $q$-character of the Kirillov-Reshetikhin
module $W_{k+1,q^r}$ is given by
\begin{align*}
 \chi_q(W_{k+1,q^r})&= Y_{q^r} Y_{q^{r+2}} \cdots Y_{q^{r+2k}}
\left(1+A_{q^{r+2k+1}}^{-1}+A_{q^{r+2k-1}}^{-1}A_{q^{r+2k+1}}^{-1}
+ \right.
\\
&\left.\ \ \ +\ \cdots +A_{q^{r+1}}^{-1}A_{q^{r+3}}^{-1}\cdots A_{q^{r+2k+1}}^{-1}\right),
\end{align*}
in agreement with Eq.~(\ref{eqKR}).

\subsection{Algorithms}

Let $m\in\M_+$. We say that the simple module $L(m)$ is minuscule 
if $m$ is the only dominant monomial of $\chi_q(L(m))$. 
(In \cite{N3} these modules are called special.)
There exists 
an algorithm due to Frenkel and Mukhin \cite{FM} which
attaches to any $m\in\M_+$ a polynomial $\FM(m)\in\Y$, and
in case $L(m)$ is minuscule it is proved that $\FM(m) = \chi_q(L(m))$.
Moreover, all fundamental modules $L(Y_{i,a})$ are minuscule,
so this algorithm allows to calculate their
$q$-characters.  
It was proved in \cite{N2} that Kirillov-Reshetikhin modules are also
minuscule.
But there also exist simple modules for which the Frenkel-Mukhin
algorithm fails. For example in type~$A_2$, 
$
\chi_q(L(Y_{1,1}^2Y_{2,q^3})) \not = \FM(Y_{1,1}^2Y_{2,q^3}),
$
\cite[Example 5.6]{HL}. (For an earlier example in type $C_3$
see \cite{NN}.)

In \cite{N3}, Nakajima has introduced a $t$-analogue $\chi_{q,t}$
of the $q$-character $\chi_q$. This is obtained by keeping the $t$-grading 
in the graded multiplicities $a_{V,0;W}(t)$ of Theorem~\ref{theo_sim}.
Imitating the Kazhdan-Lusztig algorithm for calculating the intersection 
cohomology of a Schubert variety, he has described an algorithm for
computing the $(q,t)$-character of an arbitrary simple module in terms
of the $(q,t)$-characters of the fundamental modules.
The $(q,t)$-characters of the fundamental modules can in turn
be obtained using a $t$-version of the Frenkel-Mukhin algorithm.
We therefore have, in principle, a way of calculating $\chi_q(L(m))$
for every $m\in\M_+$.


\section{Tensor structure}\label{sect_4}

In the category of finite-dimensional $U_q(\gg)$-modules, 
tensor products of irreducible modules are almost never 
irreducible. This is in sharp contrast with what happens 
for tensor products of finite-dimensional $U_q(L\gg)$-modules.
Indeed, if $M$ and $N$ are simple objects of $\CC$, the tensor
product $M\otimes N(a)$ is simple for all but a finite
number of $a\in\C^*$. (Here $N(a)$ is the image of 
$N$ under the auto-equivalence $\tau_a^*$ of \S\ref{subsect_2.1}.)  
Hence many tensor products of simple $U_q(L\gg)$-modules are
simple, or equivalently, many simple modules can be factored
as tensor products of smaller simple modules.
The following questions are therefore natural:
\begin{itemize}
 \item[(i)] what are the \emph{prime} simple modules, \ie 
the simple modules which have no factorization
as a tensor product of smaller modules ?
 \item[(ii)] which tensor products of prime simples are
simple ? 
\end{itemize}
We have seen in \S\ref{subsect_2.4} that these questions have
a simple answer when $\gg=\Sl_2$, namely, the prime simples
are the Kirillov-Reshetikhin modules, and a tensor product
of Kirillov-Reshetikhin modules is simple if and only if
the corresponding $q$-segments are pairwise in general position.
In this third lecture, we will report on some recent progress
in trying to extend these results to an arbitrary simply-laced~$\gg$.

\subsection{The cluster algebra $\A_\ell$}
\begin{figure}[t]
\[
\def\objectstyle{\scriptstyle}
\def\lablestyle{\scriptstyle}
\xymatrix@-1.0pc{
&&\ar[ld] {(2,5)} \ar[rd]&&
\\
&{(1,4)}\ar[rd]&
&\ar[ld] {(3,4)} 
\\
&&\ar[ld] (2,3)\ar[uu] \ar[rd]&&
\\
&(1,2)\ar[uu] \ar[rd] &&\ar[ld] (3,2)\ar[uu]
\\
&&\ar[ld] (2,1)\ar[uu] \ar[rd]&&
\\
&(1,0)\ar[uu] && (3,0)\ar[uu] 
\\
}
\]
\caption{\label{Fig4} {\it The quiver $\Ga_2$ in type $A_3$.}}
\end{figure}
We will assume the reader has some familiarity with cluster
algebras. Nice introductions to this theory with pointers to the
literature have been written by Zelevinsky \cite{Z} and Fomin \cite{F}.
All the necessary material for understanding this lecture can
also be found in \cite[\S2]{Le}.

For $\ell\in\N$, we define a new quiver $\Ga_\ell$.
Put 
\[
\hI_0(\ell)=\{(i,\xi_i+2k)\mid i\in I,\ 0\le k \le \ell\}.
\]
The graph $\Ga_\ell$ is obtained by taking the full subgraph of $\Z Q$ with
vertex set $\hI_0(\ell)$, and by adding to it new vertical up-arrows 
corresponding to the natural translation $(i,r) \mapsto (i,r+2)$.
For example, if $\gg$ has type $A_3$ and $I_0 = \{1,3\}$, the quiver $\Ga_2$ is
shown in Figure~\ref{Fig4}.

\begin{table}[t]
\begin{center}
\begin{tabular}
{|c|c|c|}
\hline
Type of $\gg$ &  $\ell$ & Type of $\A_\ell$\\
\hline
$A_1$ & $\ell$ & $A_\ell$ \\
\hline
$X_n$ & $1$ & $X_n$ \\
\hline
$A_2$ & $2$ & $D_4$ \\
$A_2$ & $3$ & $E_6$\\
$A_2$ & $4$ & $E_8$ \\
\hline
$A_3$ & $2$ & $E_6$ \\
\hline
$A_4$ & $2$ & $E_8$ \\
\hline
\end{tabular}
\end{center}
\caption{\small \it Algebras ${\cal A}_\ell$ of finite cluster type.
\label{table1}}
\end{table}

Let $\z = \{z_{(i,r)} \mid (i,r)\in \hI_0(\ell)\}$ be a set of indeterminates
corresponding to the vertices of $\Ga_\ell$,
and consider the seed $(\z , \Ga_\ell)$ in which the variables
$z_{(i,\xi_i)}\ (i\in I)$ are frozen.
This is the initial seed of a cluster algebra $\A_\ell \subset \Q(\z)$.
It follows easily from \cite{FZ2} that $\A_\ell$ has in general
infinitely many cluster variables.
The exceptional pairs $(\gg,\ell)$
for which $\A_\ell$ has finite cluster type are listed in Table~\ref{table1}.

\subsection{Conjectural relation between $\A_\ell$ and $\CC_\ell$}

Recall the subcategory $\CC_\ell$ from \S\ref{subsect_2.5}, and its Grothendieck
ring $R_\ell$.
We say that a simple object $S$ of $\CC_\ell$ is \emph{real} if $S\otimes S$
is simple. 
 
\begin{conjecture}[Hernandez-Leclerc \cite{HL}]\label{conjec}
The assignment 
\[
z_{(i,\xi_i+2k)} \mapsto \left[W^{(i)}_{\ell+1-k,\,q^{\xi_i+2k}}\right]
\]
extends to a ring isomorphism $\iota_\ell: \A_{\ell} \to R_\ell$.
The map $\iota_\ell$ induces a bijection between
cluster monomials and classes of real simple objects of $\CC_\ell$, 
and between cluster variables and
classes of real prime simple objects of $\CC_\ell$.
\end{conjecture}

Note that in \cite{HL} we have chosen a different initial seed
for defining $\A_\ell$, so the Kirillov-Reshetikhin
modules assigned to the initial cluster variables have different
spectral parameters\footnote{We take this opportunity to correct
a typo in \cite{Le}: in the statement of Conjecture 9.1, 
one should replace $W^{(i)}_{k,\,q^{\xi_i+2(\ell+1-k)}}$ by $W^{(i)}_{k,\,q^{\xi_i}}$.}.
 
For $\gg=\Sl_2$, Conjecture~\ref{conjec} holds, as it is just a reformulation
of the classical results of \S\ref{subsect_2.4}.
If true in general, Conjecture~\ref{conjec} will give a combinatorial description in terms of
cluster algebras of the prime tensor factorization of every real
simple module of~$\CC$.

Note that, by definition, the square of a cluster monomial is again a cluster monomial.
This explains why cluster monomials can only correspond to real simple modules.
For $\gg=\Sl_2$, all simple $U_q(L\gg)$-modules are real.
However for $\gg\not = \Sl_2$ there exist \emph{imaginary} simple $U_q(L\gg)$-modules
(\ie simple modules whose tensor square is not simple), as shown
in \cite{Le1}.
This is consistent with the expectation that 
a cluster algebra with infinitely many cluster variables
is not spanned by its set of cluster monomials.

We arrived at Conjecture~\ref{conjec} by noting that the $T$-system
equations satisfied by Kirillov-Reshetikhin modules (see \S\ref{subsect_2.3})
are of the same form as the exchange relations of a cluster algebra.
This was inspired by the seminal work \cite{FZ0},
in which cluster algebra combinatorics is used to prove Zamolodchikov's 
periodicity conjecture for $Y$-systems attached to Dynkin diagrams.

\subsection{The case $\ell = 1$}

Our main evidence for Conjecture~\ref{conjec} is the following:
\begin{theorem}[\cite{HL,N3}]\label{th_l_1}
Conjecture~\ref{conjec} holds for $\gg$ of type $A, D, E$ and $\ell = 1$.
In this case, all simple modules are real. 
\end{theorem}
This was first proved in \cite{HL} for type $A$ and $D_4$ by combinatorial
and represen\-ta\-tion-theoretic methods, and soon after, by Nakajima \cite{N3} 
in the general case, by using the geometric description of the irreducible
$q$-characters explained in \S\ref{sect_3}.
These two different proofs will be explained in \S\ref{proofHL} and \S\ref{proofN}.

Let us illustrate Theorem~\ref{th_l_1} for $\gg = \Sl_4$.
As the cluster algebra $\A_1$ has finite cluster type $A_3$,
the cluster variables, and therefore the non frozen prime simple modules,
are in bijection with the almost positive roots of $A_3$ \cite{FZ2}.
Of course, there can be several such bijections. The bijection
chosen in \cite{HL} is as follows:
\begin{align*}
&S(-\al_1)= L(Y_{1,q^2}),\ 
S(-\al_2) = L(Y_{2,q}),\
S(-\al_3) = L(Y_{3,q^2}),\\
&S(\al_1) = L(Y_{1,q^0}),\ \ \
S(\al_2) = L(Y_{2,q^3}),\ \ \
S(\al_3) = L(Y_{3,q^0}),\\ 
&S(\al_1+\al_2) = L(Y_{1,q^0}Y_{2,q^3}),\ \
S(\al_2+\al_3) = L(Y_{2,q^3}Y_{3,q^0}),\\
&S(\al_1+\al_2+\al_3) = L(Y_{1,q^0}Y_{2,q^3}Y_{3,q^0}).  
\end{align*}
Note that the last three modules are not Kirillov-Reshetikhin modules. 
There are three more prime simples corresponding to the three frozen
variables of $\A_1$, namely
\[
F_1 = L(Y_{1,q^0}Y_{1,q^2}),\ \
F_2 = L(Y_{2,q}Y_{2,q^3}),\ \
F_3 = L(Y_{3,q^0}Y_{3,q^2}). 
\]
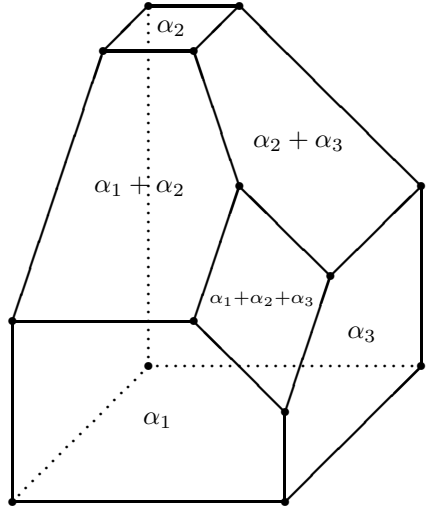
\begin{figure}[t]
\begin{center} 
\setlength{\unitlength}{1.70pt} 
\begin{picture}(90,110)(0,0) 
\thicklines 
\put(0,0){\line(1,0){60}} 
\put(0,0){\line(0,1){40}} 
\put(60,0){\line(0,1){20}} 
\put(60,0){\line(1,1){30}} 
\put(0,40){\line(1,0){40}} 
\put(0,40){\line(1,3){20}} 
\put(60,20){\line(-1,1){20}} 
\put(60,20){\line(1,3){10}} 
\put(90,30){\line(0,1){40}} 
\put(40,40){\line(1,3){10}} 
\put(70,50){\line(1,1){20}} 
\put(70,50){\line(-1,1){20}} 
\put(50,70){\line(-1,3){10}} 
\put(90,70){\line(-1,1){40}} 
\put(20,100){\line(1,0){20}} 
\put(20,100){\line(1,1){10}} 
\put(30,110){\line(1,0){20}} 
\put(40,100){\line(1,1){10}} 
 
\thinlines 
\multiput(0,0)(1.5,1.5){20}{\circle*{0.5}} 
\multiput(30,30)(2,0){30}{\circle*{0.5}} 
\multiput(30,30)(0,2){40}{\circle*{0.5}} 
 
\put(35,105){\makebox(0,0){$\alpha_2$}} 
\put(28,70){\makebox(0,0){$\alpha_1+\alpha_2$}} 
\put(63,80){\makebox(0,0){$\alpha_2+\alpha_3$}} 
\put(55,45){\makebox(0,0){$\scriptstyle \alpha_1+\alpha_2+\alpha_3$}} 
\put(32,18){\makebox(0,0){$\alpha_1$}} 
\put(77,37){\makebox(0,0){$\alpha_3$}} 
 
\put(0,0){\circle*{2}} 
\put(60,0){\circle*{2}} 
\put(60,20){\circle*{2}} 
\put(30,30){\circle*{2}} 
\put(90,30){\circle*{2}} 
\put(0,40){\circle*{2}} 
\put(40,40){\circle*{2}} 
\put(70,50){\circle*{2}} 
\put(50,70){\circle*{2}} 
\put(90,70){\circle*{2}} 
\put(20,100){\circle*{2}} 
\put(40,100){\circle*{2}} 
\put(30,110){\circle*{2}} 
\put(50,110){\circle*{2}} 
\end{picture}  
\end{center}
\caption{\label{Figassocia} {\it The associahedron of type $A_3$.}}
\end{figure}
The cluster algebra $\A_1$ has fourteen clusters, in bijection with
the vertices of the associahedron shown in Figure~\ref{Figassocia} \cite{FZ2}.
The faces of the associahedron are naturally labeled by the almost
positive roots (the rear, bottom, and leftmost faces are labeled by
$-\al_1$, $-\al_2$, and $-\al_3$, respectively). Each vertex corresponds
to the cluster consisting in its three adjacent faces:
\[
 \begin{array}{c}
\{-\al_1,\,-\al_2,\,-\al_3\},\ 
\{\al_1,\,-\al_2,\,-\al_3\},\ 
\{-\al_1,\,\al_2,\,-\al_3\},\ 
\{-\al_1,\,-\al_2,\,\al_3\}, 
\\[2mm]
\{\al_1,\,-\al_2,\,\al_3\},\,   
\{-\al_1,\,\al_2,\,\al_2+\al_3\},\,
\{-\al_1,\,\al_3,\,\al_2+\al_3\},\,
\{-\al_3,\,\al_2,\,\al_1+\al_2\},
\\[2mm]
\{-\al_3,\,\al_1,\,\al_1+\al_2\},\
\{\al_1+\al_2,\,\al_2,\,\al_2+\al_3\},\
\{\al_1,\,\al_3,\,\al_1+\al_2+\al_3\},
\\[2mm]
\{\al_1,\,\al_1+\al_2,\,\al_1+\al_2+\al_3\},\ 
\{\al_3,\,\al_2+\al_3,\,\al_1+\al_2+\al_3\}, 
\\[2mm]
\{\al_1+\al_2,\,\al_2+\al_3,\,\al_1+\al_2+\al_3\}.
\end{array}
\]
The simple modules of $\CC_1$ are exactly
all tensor products of the form
\[
S(\be_1)^{\otimes k_1}\otimes S(\be_2)^{\otimes k_2}\otimes S(\be_3)^{\otimes k_3}\otimes 
F_1^{\otimes l_1}\otimes F_2^{\otimes l_2}\otimes F_3^{\otimes l_3},
\ (k_1,k_2,k_3,l_1,l_2,l_3)\in\N^6,
\]
in which $\{\be_1,\,\be_2,\,\be_3\}$ runs over the 14 clusters listed above. 

Note that by Gabriel's theorem, positive roots are in one-to-one correspondence
with indecomposable representations of $Q$. By inspection
of the above list of clusters, one can check that two roots belong to a common cluster
if and only if the corresponding representations of $Q$ have no extension
between them. This is true in general, as will be explained below (see
Corollary~\ref{thgeom}, and the end of \S\ref{proofN}).

\subsection{Proof of Theorem~\ref{th_l_1}: approach of \cite{HL}}\label{proofHL}	

Write for short $z_i = z_{(i,\xi_i+2)}$ for the cluster variables of
the initial cluster $\z$.
We know that $R_1$ is the polynomial ring in the classes of the
fundamental modules $L(Y_{i,q^{\xi_i}})$ and $L(Y_{i,q^{\xi_i+2}})$.
On the other hand, it is not difficult to show that, because
of the presence of the frozen variables, $\A_1$
is the polynomial ring in the variables $z_i$, together with the 
variables $z_i'$ of the cluster $\z'$ obtained from $\z$
by applying the product of mutations
\[
\z' = \left(\prod_{i\in I_0}\mu_i\right)\left( \prod_{i\in I_1} \mu_i\right)\z.
\]
Therefore, the assignment 
\[
 z_i\mapsto [L(Y_{i,q^{\xi_i+2}})],\quad z'_i \mapsto [L(Y_{i,q^{\xi_i}})],\qquad (i\in I),
\]
extends to a ring isomorphism $\iota$ from $\A_1$ to $R_1$. 

To calculate the images under $\iota$ of the remaining cluster variables,
we use the fact \cite{FZ4} that every cluster variable is entirely determined by
its $F$-polynomial (and its $g$-vector) with respect to the reference cluster $\z$. 
For an almost positive root $\be=\sum_i b_i\al_i$, let $\z[\be]$ be the cluster variable 
whose cluster expansion with respect to $\z$ has denominator $\prod_i z_i^{b_i}$.
In particular $z_i = \z[-\al_i]$. Denote by $F_\be$ the $F$-polynomial 
of $\z[\be]$. By convention $F_{-\al_i}=1$.

On the other hand, the $q$-character of an object $M$ of $\CC_1$ is uniquely
determined by its truncation obtained by specializing $A_{i,q^{k}}^{-1}$ to $0$
for $k>\xi_i+1$. 
The truncated $q$-character of a simple object $L(m)$ of $\CC_1$
is of the form
\begin{equation}\label{qrest}
\chi_q(L(m))_{\le 2} = m\, P_m(v_1,\ldots,v_n) 
\end{equation}
where $P$ is a polynomial in the variables $v_i := A_{i,q^{\xi_i+1}}^{-1}\ (i\in I)$
with constant term~1. 
Moreover, the map $\tau: [L(m)]\mapsto \chi_q(L(m))_{\le 2}$ is an injective ring
homomorphism from $R_1$ to its image in $\Y$.
The injectivity comes from the fact that the truncated $q$-character of a module of $\CC_1$
already contains all its dominant monomials.
It is thus enough to determine the images of the cluster variables of $\A_1$
under $\iota':=\tau\iota$.

Let $s_i\ (i\in I)$ be the Coxeter generators of the Weyl group.
In \cite{HL}, it is proved that for $\be>0$,
\begin{equation}\label{imF}
\iota'(\z[\be]) = Y^{\al} F_\beta(v_1,\ldots,v_n),  
\end{equation}
where 
\begin{equation}\label{albe}
\al = \sum_ia_i\al_i := \left(\prod_{i\in I_1} s_i\right)\be
\end{equation}
and 
\begin{equation}\label{Yal}
Y^\al = 
\left\{
\begin{array}{ll}
\ds\prod_{i\in I} Y_{i,3\xi_i}^{a_i}&\quad \mbox{if}\quad \al>0,
\\[5mm]
Y_{i,2-\xi_i}&\quad \mbox{if} \quad \al=-\al_i.
\end{array}
\right.
\end{equation}
Thus, setting $m=Y^\al$ and comparing (\ref{imF}) with (\ref{qrest}), we see
that an important step in proving Theorem~\ref{th_l_1} is to show that
the two polynomials $P_m$ and $F_\be$ coincide. This last statement is
verified in \cite{HL} for every root $\be$ in types $A_n$ and $D_n$,
and for every multiplicity-free root $\be$ in type $E_n$.
The proof uses the Frenkel-Mukhin algorithm for evaluating $P_m$,
and the combinatorial description of the Fibonacci polynomials
of Fomin and Zelevinsky \cite{FZ1} for evaluating $F_\be$.  
Thus, except for these missing roots in type $E_n$, this shows that
all cluster variables of $\A_1$ are mapped by $\iota$ to the classes of some
simple modules in $R_1$. 

The second main step is the following tensor product theorem,
proved for all types $A_n$, $D_n$, $E_n$. 
Let $S_1,\ldots,S_k$ be simple modules of $\CC_1$, and suppose that
for every $1\le i<j\le k$ the tensor product $S_i\otimes S_j$
is simple. Then it is shown \cite[Th. 8.1]{HL} that $S_1\otimes\cdots\otimes S_k$
is simple\footnote{This theorem was later extended by Hernandez \cite{H3}
to the whole category $\CC$.}.
Thus, to show that the image of a cluster monomial by $\iota$ is the
class of a simple module, it is enough to prove it when the monomial
is the product of \emph{two} cluster variables.

Finally, the third step consists in proving that if $z[\be]$ and
$z[\ga]$ are two compatible cluster variables of $\A_1$, that is, if
$z[\be]z[\ga]$ is a cluster monomial,
then the tensor product of the corresponding simple modules of $\CC_1$
is simple.
Since, for a given $\gg$, there are only finitely many cluster
variables in $\A_1$, and so finitely many compatible pairs, this
is in principle only a ``finite check''. Unfortunately it is not easy in general
to decide if a product of (truncated) irreducible $q$-characters
is simple, and in \cite{HL} this was only settled completely in types
$A_n$ and~$D_4$. 

Although this (partial) proof is combinatorial and representation-theoretic,
it has an interesting geometric consequence. Indeed, it shows that the
truncated $q$-characters of the prime simple objects of $\CC_1$ 
coincide, after dividing out the highest {\em l}-weight monomial, with the $F$-polynomials
of the cluster variables of $\A_1$. But the $F$-polynomials have a 
geometric description due to Fu and Keller \cite{FK} in terms of
quiver Grassmannians, inspired from a similar formula 
of Caldero and Chapoton for cluster expansions of cluster variables \cite{CC}. 
Therefore we get the following geometric description of the truncated
$q$-characters.

Let $M[\be]$ be the indecomposable representation of the Dynkin quiver $Q$
attached to a positive root $\be$, and denote by
$\Gr_\nu(M[\be])$ the quiver Grassmannian of subrepresentations of $M[\be]$
with dimension vector $\nu$.
Let $\al$ and $Y^\al$ be related to~$\be$ as in Eq.(\ref{albe}), (\ref{Yal}).
Finally, recall the notation $v_i := A_{i,q^{\xi_i+1}}^{-1}$.
\begin{corollary}[\cite{HL}]\label{thgeom}
Conjecture~\ref{conjec} for $\CC_1$ implies that
\begin{equation}\label{quivgrass}
\chi_q(L(Y^\al))_{\le 2}\ = 
Y^\al \sum_\nu \chi(\Gr_\nu(M[\be])) \,v_1^{\nu_1}\cdots v_n^{\nu_n}.
\end{equation}
More generally, we have a similar truncated $q$-character formula for
every simple module of $\CC_1$, in which the indecomposable
representation $M[\be]$ of the right-hand side is replaced by a generic representation of $Q$,
\ie a representation without self-extension. 
\end{corollary}
Corollary~\ref{thgeom} should be compared to Eq.(\ref{eqstandard}) and (\ref{eqfund})
for $q$-characters of standard modules. What is remarkable here is that
we obtain a similar formula for \emph{simple} modules of $\CC_1$: for all these
modules, we do not need to use the decomposition theorem for perverse sheaves,
as was done in \S\ref{perverse}.

\subsection{Proof of Theorem~\ref{th_l_1}: approach of \cite{N5}}\label{proofN}
\begin{figure}[t]
\[
\xymatrix@-1.0pc{
&\qquad\qquad\qquad&&
\\
&\ar[d]^{\al_2(3)}W_2(3)&
\\
\ar[d]^{\al_1(2)} W_1(2)&\ar[ld]^{B_{21}(2)} V_2(2)\ar[d]^{\be_2(2)} \ar[rd]^{B_{23}(2)} & \ar[d]^{\al_3(2)}W_3(2)&
\\
V_1(1)\ar[d]^{\be_1(1)}  &W_2(1) & V_3(1)\ar[d]^{\be_3(1)}
\\
W_1(0) && W_3(0)
}
\]
\caption{\label{Fig5} {\it The graded spaces $W$ and $V$ associated with a simple object of $\CC_1$ in type $A_3$.}}
\end{figure}
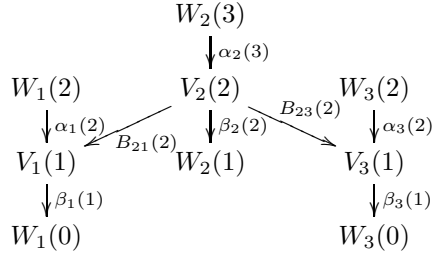
In \cite{N5}, Nakajima reverses the logic of \cite{HL}, and 
first proves the formula of
Corollary~\ref{thgeom} for all simple modules 
of $\CC_1$ and for all Dynkin types, by means of his
description in terms of perverse sheaves.
This is made possible because of the following simple description
of the quiver varieties $\fM^\bullet(V,W)$ and $\fM_0^\bullet(W)$
when $W$ corresponds to the highest {\em l}-weight of a simple 
object of $\CC_1$, and the monomial $Y^WA^V$ contributes to its truncated
$q$-character.

One first notes that $L(W)$ is in $\CC_1$ if and only if the 
$\hI_0$-graded space $W$ satisfies:
\begin{equation}
 W_i(r) \not = 0 \quad\mbox{ only if }\quad r\in\{\xi_i, \xi_i+2\}.
\end{equation}
Moreover, if $Y^WA^V$ appears in 
$\chi_q(L(W))_{\le 2}$ then 
\begin{equation}\label{eqV}
 V_i(r) \not = 0 \quad\mbox{ only if }\quad r=\xi_i+1.
\end{equation}
Thus $W$ and $V$ are supported on a zig-zag strip of height 2,
as shown in Figure~\ref{Fig5}.
Therefore the ADHM equations are trivially satisfied in this case,
and we have $M^\bullet(V,W)=\L^\bullet(V,W)$.

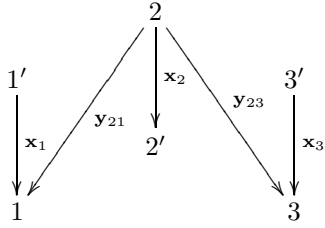
\begin{figure}[t]
\[
\xymatrix@-1.0pc{
&\qquad\qquad\qquad&&
\\
&\ar[lddd]^{\y_{21}}2\ar[dd]^{\x_2}\ar[rddd]^{\y_{23}}&
\\
1'\ar[dd]^{\x_1}& & 3'\ar[dd]^{\x_3}&
\\
&2' & 
\\
1&& 3
}
\]
\caption{\label{Fig6} {\it The decorated quiver $\tQ$ in type $A_3$.}}
\end{figure}

Since every dominant monomial appears in the truncated $q$-character,
and since $\fM_0^\bullet(W)$ is equal to $\fM_0^\bullet(V,W)$
for some {\em l}-dominant pair $(V,W)$, we see that 
$\fM_0^\bullet(W) = M^\bullet(V,W)\sslash G_V$ for some $V$
satisfying (\ref{eqV}).
Define the following $G_V$-invariant maps on $M^\bullet(V,W)$:
\begin{align}
\x_i &= \be_i(\xi_i+1)\al_i(\xi_i+2),\quad (i\in I),\label{eq28}
\\
\y_{ij} &=  \be_j(1)B_{ij}(2)\al_i(3),\quad\quad (i\in I_1,\ j\in I_0,\ a_{ij}=-1).\label{eq29}
\end{align}
The data $(\x_i,\y_{ij})$ amount to a representation of a decorated quiver
$\tQ$ obtained by attaching to every vertex $i$ of $Q$ a new vertex $i'$
and an arrow $i'\to i$ (\resp $i\to i'$) if $i\in I_0$ (\resp $i\in I_1$)
(see Figure~\ref{Fig6}). Let
\[
 E_W = \bigoplus_{i\in I} \Hom(W_i(\xi_i + 2),W_i(\xi_i)) \oplus
\bigoplus_{i\in I_1,\, j\in I_0,\, a_{ij}=-1} \Hom(W_i(3),W_j(0))
\]
be the space of representations of $\tQ$ based on $W$.
Nakajima shows that the map $(B,\al,\be) \mapsto (\x_i, \y_{ij})$ 
induces an isomorphism from $\fM_0^\bullet(W)$ to $E_W$.
Hence, for $L(W)$ in $\CC_1$ the affine variety $\fM_0^\bullet(W)$
is isomorphic to a vector space.

Put $\nu_i = \dim V_i$ and $\nu = (\nu_i) \in \N^I$.
Let $\F(\nu,W)$ be the variety of $n$-tuples $X=(X_i)$
of subspaces of $W$ satisfying $\dim X_i = \nu_i$ and
\[
X_i \subseteq W_i(0) \quad (i\in I_0),
\qquad
X_i \subseteq W_i(1) \oplus \bigoplus_{j:\, a_{ij}=-1} X_j \quad (i\in I_1). 
\]
Define $\tF(\nu,W)$ as the closed subvariety of $E_W \times \F(\nu,W)$ consisting
of all elements $((\x_i,\,\y_{ij}),\, X)$ such that
\begin{equation}\label{fiber}
\im \x_i \subseteq X_i\ (i\in I_0),\qquad 
\im\left(\x_i\oplus \bigoplus_{j:\, a_{ij}=-1} \y_{ij}\right) \subseteq X_i
\ (i\in I_1).
\end{equation}
Nakajima shows that $(B,\al,\be)\in M^\bullet(V,W)$ is stable if and only if
all maps $\be_i(1)\ (i\in I_0)$ and $\si_i(2) := \be_i(2) + 
\sum_{j:\,a_{ij}=-1} B_{ij}(2)\ (i\in I_1)$ are injective.
Clearly the collection $X$ of spaces 
\begin{align}
X_i& = \be_i(1)(V_i(1))\ (i\in I_0),\label{eq30}
\\ 
X_i& = \si_i(2)(V_i(2))\ (i\in I_1), \label{eq31}
\end{align}
is $G_V$-invariant, and $\dim X_i = \nu_i$ if $(B,\al,\be)$ is stable.
Therefore, the map $(B,\al,\be)\mapsto ((\x_i,\,\y_{ij}),\, X)$ defined by
(\ref{eq28}) (\ref{eq29}) (\ref{eq30}) (\ref{eq31}) induces a map from
$\fM^\bullet(V,W)$ to $\tF(\nu,W)$, and Nakajima shows that this is
an isomorphism \cite[Proposition 4.6]{N5}.
Moreover, when $\fM^\bullet(V,W)$ and $\fM_0^\bullet(W)$ are realized
as $\tF(\nu,W)$ and $E_W$, respectively, then the projective morphism
$\pi_V: \fM^\bullet(V,W) \to \fM_0^\bullet(W)$ becomes the first projection.
(Compare this description with the prototypical example of \S\ref{subsect_3.5}.1.) 

By Theorem~\ref{theo_sim}, to calculate the truncated $q$-character of
a simple module of~$\CC_1$, one must now compute the multiplicity $a_{V,0;W}(1)$ of
the skyscraper sheaf $IC_W(0)$ in the expansion of 
$[(\pi_V)_!(1_{\fM^\bullet(V,W)})]$  
on the basis $\{[IC_W(V')]\}$.
Since $\fM_0^\bullet(W) \simeq E_W$ is a vector space, one
can use for that a Fourier transform. 
Let $E_W^*$ denote the dual space, and let $\psi$ be the
Fourier-Sato-Deligne functor from the derived category 
$D(E_W)$ to $D(E_W^*)$. The functor $\psi$ maps every 
simple perverse sheaf $IC_W(V)$ on $E_W$ to a simple
perverse sheaf on $E_W^*$. In particular, the image of the
skyscraper sheaf is 
\[
 \psi(IC_W(0)) = 1_{E_W^*}[\dim E_W],
\]
the constant sheaf on $E_W^*$, with degree shifted by $\dim E_W$.

We can regard the product $E_W \times \F(\nu,W)$ as a
trivial vector bundle on $\F(\nu,W)$ with fiber $E_W$.
By (\ref{fiber}), the fibers of the restriction of the second projection to 
$\tF(\nu,W)$ are vector spaces of constant dimension, and
$\tF(\nu,W)$ can be seen as a subbundle of $E_W \times \F(\nu,W)$.
Denote by $\tF(\nu,W)^\perp$ the annihilator of
$\tF(\nu,W)$ in the dual trivial bundle $E_W^* \times \F(\nu,W)$.

We also have a Fourier-Sato-Deligne functor $\psi'$ from the 
derived category of the trivial bundle $E_W \times \F(\nu,W)$
to that of $E_W^* \times \F(\nu,W)$.
It satisfies
\[
 \psi'\left(1_{\tF(\nu,W)}[\dim \tF(\nu,W)]\right) = 1_{\tF(\nu,W)^\perp}[\dim \tF(\nu,W)^\perp].
\]
Moreover, denoting by $\pi: \tF(\nu,W)\to E_W$ and $\pi^\perp: \tF(\nu,W)^\perp\to E_W^*$ the bundle maps,
we have the commutation relation
\[
\pi_!^\perp\circ \psi' = \psi \circ \pi_!. 
\]
It follows that the required (ungraded) multiplicity $a_{V,0;W}(1)$ is equal
to the multiplicity of the constant sheaf $1_{E_W^*}$ in the expansion of
$\pi_!^\perp (1_{\tF(\nu,W)^\perp})$ in terms of the $\{\psi(IC_W(V))\}$.
The advantage of this Fourier transformation is that we can now evaluate this
new multiplicity by looking at the stalk of $\pi_!^\perp (1_{\tF(\nu,W)^\perp})$
over a generic point of $E_W^*$.

At this point we remark that we can without loss of generality assume
that $W_i(2-\xi_i)=0$ for every $i\in I$. In other words, we suppose
that $E_W$ is a space of representations of the quiver $Q$ \emph{without decoration}.
(One can easily reduce the general case to this one by factoring out from
$L(W)$ a tensor product of frozen Kirillov-Reshetikhin modules
$L(Y_{i,\xi_i} Y_{i,\xi_i+2})$, as in \cite[\S9.2]{HL} or \cite[\S6.3]{N5}.)  
So $E_W^*$ is a space of representations of the quiver $Q^*$ obtained from $Q$
by changing the orientation of every arrow.
Let 
\[
 G_W := \prod_{i\in I} GL(W_i(3\xi_i)). 
\]
Since $Q^*$ is a Dynkin quiver, $E_W^*$ has an open dense $G_W$-orbit corresponding
to the generic representation of dimension vector $(\dim W_i(3\xi_i))$, and
all other $G_W$-orbits have strictly smaller dimension. 
Now, all the simple perverse sheaves $\psi(IC_W(V))$ are $G_W$-equivariant,
hence they are supported on a union of $G_W$-orbits, 
so the only one having a nonzero stalk over a generic point of $E_W^*$
is $\psi(IC_W(0)) = 1_{E_W^*}[\dim E_W]$. Therefore, by definition of the
pushdown functor~$\pi_!^\perp$, the multiplicity 
$a_{V,0;W}(1)$ is nothing else than the dimension of the total cohomology
of a generic fiber of $\pi^\perp$.

It remains to describe this generic fiber.  
Because of our simplifying assumption, a point of $E_W$ is now just
a collection of maps $\y_{ij} \in \Hom(W_i(3),W_j(0))\ (a_{ij}=-1)$,
and a point in $\F(\nu,W)$ is a collection of subspaces $X=(X_i)$
of $W$ such that 
\[
X_i \subseteq W_i(0)\ \ (i\in I_0), \qquad
X_i \subseteq \bigoplus_{j: a_{ij}=-1} X_j \ \ (i\in I_1).
\]
The pair $((\y_{ij}),X)$ belongs to $\tF(\nu,W)$ if and only if
$\im(\oplus_j \y_{ij}) \subseteq X_i$ for all $i\in I_1$.
Clearly, the annihilator $\tF(\nu,W)^\perp$ consists of pairs
$((\y_{ij}^*),X)$ in $E_W^*\times\F(\nu,W)$ such that 
$X_i \subseteq \Ker (\oplus_j \y_{ij}^*)$ for every $i\in I_1$.
To get a nicer description of the fibers of $\pi^\perp$ we 
consider the product $\si$ of Gelfand-Ponomarev reflection functors
at every sink $i\in I_1$ of $Q^*$. 
The functor $\si$ sends $(\y_{ij}^*)\in E_W^*$ to $(\y_{ij}^\si)\in E_{W^\si}$
defined by
\[
W_i^\si(0)=W_i(0),\ \ (i\in I_0),\qquad
W_i^\si(3)= \Ker (\oplus_j \y_{ij}^*),\ \ (i\in I_1),
\]
and, for $i\in I_1$, $\y_{ik}^\si$ is the composition of the embedding of
$\Ker (\oplus_j \y_{ij}^*)$ in $\oplus_j W_j(0)$ followed by
the projection onto $W_k(0)$.
The collection of linear maps $\y^\si=(\y_{ij}^\si)$ is a representation of
the original quiver $Q$. 
By construction, $X_i \subseteq W_i^\si(3\xi_i)$ for every $i\in I$.
Moreover, one can easily check that 
$((\y_{ij}^*),X)\in\tF(\nu,W)^\perp$ if and only if 
$\y_{ij}^\si(X_i) \subseteq X_j$ for every $i\in I_1$. 
In other words, $X$ belongs to the fiber of $\pi^\perp$ above
$(\y_{ij}^*)$ if and only if $X$ is a point of the quiver Grassmannian
$\Gr_{\nu}(\y^\si)$.

It now follows that the multiplicity $a_{V,0;W}(1)$ of the monomial
$Y^WA^V$ in $\chi_q(L(Y^W))_{\le 2}$ is the total dimension of the
cohomology of $\Gr_{\nu}(\y^\si)$ for a generic representation $\y^\si$ of $Q$
in $E_{W^\si}$. Note that the product of reflection functors $\si$ categorifies
the product $\prod_{i\in I_1} s_i$ in the Weyl group, so
if we denote by $\be$ the graded dimension 
of $W^\si$, and if we assume that $\be$ is a positive root, then
the graded dimension $\al$ of $W$ is related to $\be$ by (\ref{albe}),
in perfect agreement with (\ref{quivgrass}).
Moreover, Nakajima explains that the vanishing of the odd cohomology of 
$\fL^\bullet(V,W)$ implies that this generic fiber has no odd cohomology,
therefore $a_{V,0;W}(1)$ is also equal to the Euler characteristic of
the quiver Grassmannian $\Gr_{\nu}(\y^\si)$.
Thus, Corollary~\ref{thgeom} follows in full generality.

After this $q$-character formula is established, Nakajima proceeds to show that the tensor product factorization  
of the simple modules $L(W)$ of $\CC_1$ is given by the canonical direct sum 
decomposition of the corresponding generic quiver representation $\y$
of $E_W$ into indecomposable summands. The proof uses the geometric realization
given by Varagnolo
and Vasserot \cite{VV}
of the $t$-deformed product of $(q,t)$-characters 
in terms of convolution of perverse sheaves. 

Finally, to relate the $q$-character formula with cluster algebras, Nakajima
makes use of the cluster category of Buan, Marsh, Reineke, Reiten and Todorov \cite{BMRRT},
and of the Caldero-Chapoton formula for cluster variables \cite{CC}. 

It is worth noting that Nakajima's approach is more general: most
of his results work for the quantum affinization $U_q(L\gg)$ of a
possibly infinite-dimensional symmetric Kac-Moody algebra $\gg$. 
This yields some important positivity results for all cluster algebras
attached to an arbitrary bipartite quiver. However, when $\gg$ is infinite-dimensional
$U_q(L\gg)$ is no longer a Hopf algebra, and the
meaning of the multiplicative structure of the Grothendieck group is less clear
(see \cite{H4}).

\subsection{The case $\ell>1$}

If $\gg = \Sl_2$, Conjecture~\ref{conjec} holds for every $\ell$.
Otherwise, Conjecture~\ref{conjec} has only been proved for $\gg = \Sl_3$
and $\ell = 2$ \cite[\S13]{HL}.   
In that small rank case, $\A_2$ still has finite cluster type $D_4$
(see Table~\ref{table1}), 
and this implies that $\CC_2$ has only real objects.
There are 18 explicit prime simple objects with respective dimensions
\[
3,\ 3,\ 3,\ 3,\ 3,\ 3,\ 6,\ 6,\ 6,\ 6,\ 8,\ 8,\ 8,\ 10,\ 10,\ 15,\ 15,\ 35,  
\]
and 50 factorization patterns (corresponding to the 50 vertices of
a generalized associahedron of type $D_4$ \cite{FZ2}).
Our proof in this case is quite indirect and uses a lot of
ingredients: the quantum affine Schur-Weyl duality,
Ariki's theorem for type $A$ affine Hecke algebras \cite{A}, the
coincidence of Lusztig's dual canonical and dual semicanonical bases
of $\C[N]$ in type $A_4$ \cite{GLS1}, and the results of \cite{GLS2}
on cluster algebras and dual semicanonical bases.
This proof could be extended to $\gg = \Sl_n$ and every $\ell$ if
the general conjecture 
of \cite{GLS2} about the relationship between Lusztig's dual canonical 
and dual semicanonical bases was established.




\end{document}